\documentclass[12pt]{article}
\usepackage{amsmath,amsfonts,amsthm}

\newtheorem{thm}{Theorem}[section]
\newtheorem{lem}[thm]{Lemma}
\newtheorem{cor}[thm]{Corollary}
\newtheorem{prp}[thm]{Proposition}

\theoremstyle{definition}

\newenvironment{pfof}[1]{\begin{trivlist}\item[\hspace{\labelsep}\bf
Proof of #1 ]}{\end{trivlist}}

\newenvironment{rem}{\begin{trivlist}\item[\hspace{\labelsep}\bf
Remark ]}{\end{trivlist}}

\newcommand{\R}{\mathbb{R}}
\newcommand{\N}{\mathbb{N}}

\newcommand{\AC}{\mathrm{AC}}

\newcommand{\eps}{\varepsilon}
\renewcommand{\phi}{\varphi}

\newcommand{\dist}{\mathrm{dist}}
\newcommand{\diam}{\mathrm{diam}}
\newcommand{\Lip}{\mathrm{Lip}}

\newcommand{\iref}[1]{{\rm(\ref{#1})}}

\renewcommand{\theenumi}{\roman{enumi}}

\newcommand{\greek}[1]{\ifcase\arabic{#1}\or
$\alpha$\or$\beta$\or$\gamma$\or$\delta$\or$\varepsilon$\or$\zeta$\or
$\eta$\or$\vartheta$\or$\iota$\or$\kappa$\or$\lambda$\or$\mu$\or
$\nu$\or$\xi$\else$\dots$\fi}

\makeatletter

\newbox\bkbox

\newenvironment{in.enumerate}[1]
{\edef\clb{\@currentlabel}\begin{enumerate}\setbox\bkbox=\hbox{\clb}%
 \itemindent\wd\bkbox
 \def\atmp{#1}\def\btmp{\item}
 \ifx\atmp\btmp\def\theenumi{\clb.\arabic{\@enumctr}}\item\else
 \def\theenumi{\clb.#1{\@enumctr}}\fi}
{\end{enumerate}}

\makeatother
\def\F{{\cal L}}      
\def\E{{\cal L}}      
\def\EE{{\cal E}}     

\def\GF{\tilde\F}     
\def\GE{\tilde\E}     

\def\HH{\mathcal{H}^1_\infty} 

\def\LI{L^{\mathrm{in}}}
\def\LO{L^{\mathrm{out}}}

\begin{document}

\title{Universal singular sets in the calculus of variations}
\author{Marianna Cs\"ornyei\\ Department of Mathematics\\
University College London\\ Gower Street\\ London,  WC1E 6BT, UK
\and
Bernd Kirchheim\\ Mathematical Institute\\ University of Oxford\\\
24-29 St Giles'\\ Oxford, OX1 3LB, UK\and
Toby C.~O'Neil\thanks{Corresponding author.}\\ Department of Mathematics\\
The Open University\\ Walton Hall\\  Milton Keynes, MK7 6AA, UK \and
David Preiss\\ Mathematics Institute\\ University of Warwick\\
Coventry, CV4 7AL, UK \and   Steffen Winter\\
Friedrich-Schiller-Universit\"at Jena\\
    Mathematisches Institut\\
    D-07737 Jena\\ Germany}

\maketitle\begin{abstract}
  For regular one-dimensional variational problems, Ball and
  Nadirashvilli introduced the notion of the universal singular set of a Lagrangian
  $L$ and established
  its topological negligibility. This set
   is defined to be the set of all points in the plane through which the graph
  of some  absolutely continuous $L$-minimizer passes with infinite
  derivative.

  Motivated by Tonelli's partial regularity results, the
  question of the size of the universal singular set in measure
  naturally arises.
  Here we show that universal singular sets are characterized by being
  essentially purely unrectifiable --- that is, they intersect most
  Lipschitz curves in sets of zero length and that any compact purely unrectifiable
  set is contained within the universal singular set of some smooth Lagrangian with
  given superlinear growth. This gives examples of universal
  singular sets of  Hausdorff dimension two, filling the gap between previously known
  one-dimensional examples and Sychev's result that universal singular sets are
  Lebesgue null.

  We show that some smoothness of the Lagrangian is necessary
  for the topological size estimate, and investigate the relationship between
  growth of the Lagrangian and the existence of (pathological) rectifiable pieces in the
  universal singular set.

  We also show that Tonelli's partial regularity result is stable in that the energy of
  a `near' minimizer \(u\) over the set where it has large derivative is controlled by the how
  far \(u\) is from being a minimizer.
\end{abstract}
\tableofcontents
\section{Introduction}
In his paper~\cite{tonelli} of 1915, Tonelli gave a rigourous
treatment of the variational problem
$$ {\F} (u)=\int_a^b L(x,u(x),u'(x)) \, dx \to \min $$
over the class of absolutely continuous $u$ subject to a Dirichlet
boundary condition. It is now well known that the crucial
assumptions for attainment of the minimum are superlinear growth and
convexity of the Lagrangian $L(x,u,p)$ in $p$. We  will always
consider this the natural setting unless otherwise stated. In
addition, some smoothness or, at the least, some kind of continuity
with respect to $(x,u)$ is required, see~\cite{Serrin}. The
vectorial situation is also an active area of investigation, see for
 example~\cite{Maly}. In Section~\ref{ls} we present corresponding
results that are natural in our framework but do not lie at the core
of our later arguments.

Equally important, but perhaps more surprising, Tonelli also obtains
in~\cite{tonelli} the first partial regularity results for the
minimizer of such Lagrangians. He proves that any such minimizer has
a continuous derivative provided one allows values in the extended
real line: \(u'\in C([a,b],\R\cup\{-\infty,\infty\}) \). Hence the
singular set of such a minimizer, that is the set where it has
infinite derivative, is closed, and since $u$ is absolutely
continuous, it has measure zero.

In the converse direction,  sharper examples are given in~\cite{bm}
that are not restricted to the case of superlinear growth and Davie
shows in~\cite{davie} that, for any compact null set $E\subset\R$,
there is a smooth, convex and superlinear  Lagrangian and an
appropriate choice of boundary conditions such that any minimizer
has infinite derivatives precisely on $E$. Superlinear Lagrangians
clearly prefer bounded derivatives of $u'$ for minimizing ${\F}(u)$.
Thus the existence of singularities must be enforced by very steep
``wells'' in the $(x,u,p)$-energy landscape of $L$ --- the natural
question arises as to how these wells can be distributed.

Motivated by this, Ball and Nadirashvili introduce
in~\cite{ball_nad} the notion of the \emph{universal singular set}
of a Lagrangian \(L\): a point $(x,y)$ is in the universal singular
set of \(L\) if there is a choice of boundary conditions so that
there is a corresponding minimizer $u$ for which \(u(x)=y\) and
$|u'(x)|=+\infty$. They show that for Lagrangians of class $C^3$ the
universal singular set is a countable union of closed sets and thus
of the first Baire category. In~\cite{sychev}, Sych\"ev lowers the
smoothness assumption to $L\in C^1$ and, more importantly, shows
that the universal singular set is of zero (2-dimensional) Lebesgue
measure.

In light of these results, the question about the ``true'' size of
universal singular sets naturally arises: for example, one can ask
whether a universal singular set can have positive length or even
Hausdorff dimension larger than one. Here, we show that the key to
investigating universal singular sets is to understand their
geometric structure rather than just making size estimates. We show
in Section~\ref{ll} that for any given compact purely unrectifiable
set $E$ and any prescribed superlinear growth, there is a smooth
Lagrangian whose universal singular set contains $E$. In particular,
there are smooth Lagrangians whose universal singular sets have
Hausdorff dimension two and contain non-trivial continua. The
converse is also true in the sense that, if a given compact set
\(E\) is such that for any superlinear growth there is a
corresponding Lagrangian whose universal set contains \(E\), then
\(E\) is purely unrectifiable. On the other hand, we show that for a
prescribed superlinear growth, there is always a Lagrangian whose
universal singular set contains rectifiable pieces, although the
universal singular set is always `almost' purely unrectifiable in
the sense that it can only intersect members of a `small' class of
rectifiable curves in a set of positive length. Together these
results (Theorems \ref{thm.near-pu}--\ref{f-sigma}) imply the
results of Sych\"ev~\cite{sychev}, and Ball and
Nadirashvili~\cite{ball_nad} on the size of universal singular sets.

In fact, the result concerning the almost pure unrectifiability of
the universal singular set still holds even when we consider
Lagrangians that satisfy much weaker hypotheses than those given in
Tonelli's paper, in particular, convexity in $p$ can be dropped. In
order to have a satisfactory existence theorem for this more general
setting, we use the notion of a generalized minimizer --- an
absolutely continuous function satisfying the given boundary
conditions that is an appropriate limit of a sequence of almost
minimizers. This notion is related to, but differs slightly from,
the usual minimizers of the relaxed problem, as discussed in
Section~\ref{relax}. The difference can be made explicit by using
some instances of the so-called Lavrentiev phenomenon, where
approximation of a minimizer $u$ on $\{x:|u'(x)=+\infty\}$ by a
smooth function fails. The first example of this type was given by
Lavrentiev in~\cite{Lavrentiev}, and then simplified by~\cite{mania}
(see also~\cite{bm})
---  a survey of recent progress in understanding this behaviour is
given in~\cite{mizel}. The corresponding notion of the universal
singular set for generalized minimizers of a Lagrangian is still
`almost' purely unrectifiable, and so Sych\"ev's result that
universal singular sets have area zero holds in this broader
context. However, it is possible in this setting to construct a
continuous Lagrangian whose associated universal singular set is
residual in the plane (see Theorem~\ref{resid} of Section~\ref{ll}),
and so the result of Ball and Nadirashvili requires at least some
smoothness.

Finally, let us briefly describe the structure of the paper. In
Section~\ref{e}, apart from presenting some lower semicontinuity
results for later use,  we mainly  introduce the concepts necessary
to handle superlinear but non-convex Lagrangians. Under the weak
smoothness conditions that we assume later, our generalized
minimizers introduced in Section~\ref{gm} do not fit completely into
the usual regularization scheme --- the details are given in
Section~\ref{relax}. In Section~\ref{tr} we study Tonelli's
regularity results  for our setting (see also~\cite{CV}) and show
that for Lagrangians that are locally uniformly Lipschitz in \(y\),
regularity of minimizers is stable in the sense that the energy of
`near' minimizers on the set where they have large derivative is
controlled by how far the `near' minimizer is from being a minimizer
--- see Corollary~\ref{regul.0}. Note that Tonelli's original approach of building
auxiliary Lagrangians allows us to lower the smoothness assumptions
on the Lagrangian whereas the approach given in~\cite{bm} applies
only to Lagrangians that are $C^3$ and strictly convex in $p$ (but
not necessarily of superlinear growth). As superlinear growth is
essential for our estimates on the universal singular set, Tonelli's
approach seems more natural.

In Sections~\ref{su} and~\ref{ll} we investigate the structure of
universal singular sets --- a reader only interested in the case of
convex Lagrangians could bypass Section~\ref{e} since the results we
use from this section are well known for this situation. In
Section~\ref{ll}, after establishing the general scheme of
construction, we derive the basic examples of ``large'' universal
singular sets. In the rest of the section we complete our studies of
the interplay between the growth of $L$ and possible ``tangential''
behaviour of the universal singular  set and the link between
smoothness of $L$ and its topological size.

\subsection*{Basic notions and notations}

The set \(\AC [a,b]\) denotes the collection of all absolutely
continuous real-valued functions on the closed interval \([a,b]\),
and for \(u\in AC[a,b]\), \(u'\) denotes the derivative of \(u\).

\(\|x\|\) denotes the usual Euclidean norm of a point \(x\), and
\(\|x\|_\infty\)  denotes the sup norm of \(x\). We occasionally
also use \(\|f\|\)  to denote the \(\sup\)-norm of a function \(f\),
provided \(f\) is bounded. For \(r>0\) and \(S\subset \R^2\),
\(B(S,r)\) denotes the open \(r\)-neighbourhood of \(S\).

A set \(E\subset\R^2\) is purely unrectifiable if it meets each
Lipschitz curve \(\gamma:\R\to\R^2\) in a set of zero length:
\(\mathcal{H}^1(E\cap\{\gamma(t):t\in\R\})=0\). A set
\(E\subset\R^2\) is rectifiable if there is a countable collection
of Lipschitz curves \(\gamma_i\colon\R\to\R^2\) for which
\(\mathcal{H}^1 \left(E\setminus
\bigcup_{i=1}^\infty\gamma_i(\R)\right)=0\).

We say that a function \(L\colon\R^3\to\R\) is a \emph{Lagrangian} if:
\begin{itemize}
    \item \(L\) is bounded from below and locally bounded from above;
    \item \(L\) is Borel measurable;
    \item there is a superlinear function
      \(\omega\colon\R\to\R\) such that
    \(L(x,y,p)\geq\omega (p)\) for all \( (x,y,p)\in\R^3\).
\end{itemize}
Recall that superlinearity of $\omega$ means that
$\lim_{|p|\to\infty} \omega(p)/|p| =\infty$.

Let \(a < b\) be real numbers. For $u\in\AC[a,b]$ we let
$${\F}(u)={\F}(u; a,b):=\int_a^b L(x,u(x),u'(x))\, dx;$$
note that the integral exists thanks to the lower boundedness
and Borel measurability of $L$. Recall that
\(u\in\AC[a,b]\) is a \emph{minimizer} for the Lagrangian
\(L\) on \([a,b]\) if
$${\F}(u)= {\E}(a, u(a);b,u(b)),$$
where
$${\E}(a,A;b,B)=
\inf\left\{{\F}(v): v\in \AC [a,b], v(a)=A,\, v(b)=B\right\}.
$$
Notice that ${\E}(a,A;b,B)$ is finite; to see this it suffices to
consider the affine function that joins \( (a,A)\) to \((b,B)\). The
same argument shows that for any constant $C$, ${\E}(a,A;b,B)$ is
bounded on every bounded set of $(a,A;b,B)$ for which $a<b$ and
$|B-A|\le C|b-a|$.

To show that `almost' minimizers satisfy an approximate version of
Tonelli's partial regularity theorem, we will need to measure how
far a given function is from being a minimizer. A convenient such
measure, the \emph{excess of $u\in\AC[a,b]$ over the interval
$[a,b]$} is defined by
$$\EE(u;a,b):=\F(u;a,b)-{\E}(a,u(a);b,u(b)).$$
Of course, $\EE(u;a,b)=0$ if and only if $u$ is a minimizer. The
fact that a restriction of a minimizer to a subinterval is a
minimizer has a simple quantitative version:

\begin{lem}\label{excess}
If $a\le\alpha<\beta\le b$ then $\EE(u;\alpha,\beta)\le\EE(u;a,b)$.
\end{lem}

\begin{proof}
Assuming, as we may, that $\EE(u;a,b)<\infty$ and so also
$\F(u;a,b)<\infty$, we extend any
$v\in\AC[\alpha,\beta]$ with $v(\alpha)=u(\alpha)$ and
$v(\beta)=u(\beta)$ to $\tilde v\in \AC[a,b]$ by
$\tilde v=u$ on $[\alpha,\beta]\setminus [a,b]$.
Then
\begin{align*}
\lefteqn{\F(u;a,\alpha)+\F(v;\alpha,\beta)+\F(u;\beta,b)}\kern3cm\\
&=\F(\tilde v ;a,b)\\
&\ge  \F(u;a,b) - \EE(u;a,b)\\
&= \F(u;a,\alpha)+\F(u;\alpha,\beta)+\F(u;\beta,b)- \EE(u;a,b),
\end{align*}
giving the statement.
\end{proof}

Given a Lagrangian \(L\colon \R^3\to\R\), the
\emph{universal singular set of \(L\)}
is defined to be all points \( (x_0,y_0)\in\R^2\) for which there
are $a\leq  x_0\leq b$ with \(a<b\) and a
minimizer $u$ for \(L\) on \([a,b]\) such that
$u(x_0)=y_0$ and $u$ has infinite derivative at $x_0$.

Of course, under the above assumptions, there may well be no
minimizers for $L$ and so the notion of the universal singular set
may make little sense. In Section~\ref{gm} we define a weaker notion
of minimizers for which the existence result is nearly trivial and
for which the corresponding universal singular sets enjoy the same
smallness properties as those for ordinary minimizers.

Those of our results that do not include study of the behaviour of
minimizers when the derivative is infinite (so, in particular, our
results of Section~\ref{tr} on Tonelli's partial regularity) are
readily transferable to the vector-valued case, with essentially
identical arguments. It is not clear, however, how the universal
singular set should be defined in this situation. A natural
definition is to say that $(x_0,y_0)$ belongs to the universal
singular set if there is a minimizer $u$ such that $u(x_0)=y_0$ and
$\lim_{x\to x_0} \|u(x)-u(x_0)\|/|x-x_0|=\infty$, which is
equivalent to the standard definition in the scalar case. In this
case only one of our arguments for the scalar case extends, showing
that the graphs of absolutely continuous functions meet the
universal singular set in a set of measure zero.

\section{Generalized minimizers and regularity}\label{e}

It is obvious that under our general assumption on the Lagrangian no
existence or regularity results for minimizers can hold. One of our
goals here is to give a notion of generalized minimizers for which
the existence results hold (and are, in fact, nearly trivial), but
for which we will see later that the universal singular set is as
small as in the standard situation. We start by revising some
classical results on compactness and lower semicontinuity, then
discuss two ways of generalizing the notion of a minimizer, and
finally show that Tonelli's partial regularity results hold, under
only a mild assumption on the Lagrangian, even for these generalized
minimizers.

\subsection{Compactness and lower semicontinuity}\label{ls}

Mainly for the sake of future reference, we record here standard
arguments showing that under our assumptions on the Lagrangians, the
sets of functions with uniformly bounded Lagrangian enjoy weak
compactness properties, even if we allow changing boundary data.
Recall that we only consider Lagrangians that are Borel measurable,
lower bounded, locally bounded from above and superlinear.

\begin{prp}\label{equi}
For every real $K$ the set of functions
$$\{u\in\AC[a,b]: -K\le a<b\le K,\ {\F}(u;a,b)\le K\}$$
is uniformly equicontinuous and the set of their derivatives
$$\{u' : u\in\AC[a,b]: -K\le a<b\le K,\ {\F}(u;a,b)\le K\}$$
is equiintegrable.
\end{prp}

\begin{proof}
Suppose that $L(x,y,p)\ge\omega(p)$ where $\omega$ is superlinear
and bounded from below.
Given any $\eps>0$, let $\tau=\eps/2K$ and use that
$\omega$ is superlinear and bounded from below
to find $C$ so that $|p|\le C+ \tau \omega(p)$
for all \(p\in\R\).
If $u\in\AC[a,b]$, ${\F}(u;a,b)\le K$,
$M\subset [a,b]$ and $|M|<\eps/2C$ then
$\int_M |u'(x)|\,dx \le
\int_M (C+\tau \omega(u'(x)))\,dx\le
C|M|+\tau {\F}(u;a,b) <\eps.$
This establishes equiintegrability.
Equicontinuity follows from equiintegrability and from
$|u(\beta)-u(\alpha)|\le \int_{[\alpha,\beta]} |u'(x)|\,dx.$
\end{proof}

It is rather natural to complement this statement by a lower
semicontinuity result, for which we however need more stringent
assumptions. Recall that $u_n\in\AC[a,b]$ converge to $u\in\AC[a,b]$
weakly if and only if they converge to $u$ pointwise and the set of
$u_n'$ is equiintegrable; of course, in this situation, pointwise
convergence of the $u_n$'s is equivalent to their uniform
convergence.

In the following theorem, our extra assumptions are that the
Lagrangian is convex in \(p\) and lower semicontinuous in \(y\) for
each fixed \( (x,p)\); this should be contrasted with our results
concerning the regularity of `near' minimizers where we only require
the Lagrangian to be (locally) Lipschitz in \(y\), see
Section~\ref{tr}.

\begin{thm}\label{wlsc}
Suppose that $L(x,y,p)$ is a Lagrangian that is  convex in $p$ for
each fixed $(x,y)$, and lower semicontinuous in $y$ for each fixed
$(x,p)$. Then the map $u\in\AC[a,b]\mapsto\F(u;a,b)$ is weakly
sequentially lower semicontinuous.
\end{thm}

\begin{proof}
The theorem follows  from  Theorem~3.6 of~\cite{butt_gia} (see
also~\cite{ioffe1}) provided that we show that \(L\) is lower
semicontinuous as a function of \( (y,p)\) for fixed \(x\). This
follows easily  from the fact \(L\) is convex in \(p\), lower
semicontinuous in \(y\) and locally bounded: for fix $x\in\R$ and
suppose that $y_n\to y_0$ and $p_n\to p_0$. By the convexity of
$L(x,y_n,p)$ in~$p$, $L(x,y_n,p)\ge f_n(p)$ where
$$
f_n(p)=\begin{cases}
(p_0-p)L(x,y_n,p_0-1) + (1+(p-p_0))L(x,y_n,p_0)&\textrm{ if } p\ge p_0\\
(p-p_0)L(x,y_n,p_0+1) + (1+(p_0-p))L(x,y_n,p_0)&\textrm{ if } p\le p_0
\end{cases}
$$
Use local upper boundedness of $L$ to find a constant $C$ so that
$L(x,y_n,p_0\pm 1)\le C$ for all $n$; hence
$L(x,y_n,p)\ge f_n(p)\ge g_n(p)$ where
$$
g_n(p)=\begin{cases}
C(p_0-p) + (1+(p-p_0))L(x,y_n,p_0)&\textrm{ if } p\ge p_0\\
C(p-p_0) + (1+(p_0-p))L(x,y_n,p_0)&\textrm{ if } p\le p_0
\end{cases}
$$
Assuming, as we may, that the sequence $L(x,y_n,p_0)$ has a
(necessarily finite) limit, which is at least $L(x,y_0,p_0)$ by
lower semicontinuity for fixed $(x,p_0)$, we see that $g_n(p)$
converge uniformly on bounded sets to a continuous function $g(p)$
such that $g(p_0)\ge L(x,y_0,p_0)$. Hence $\liminf_{n\to\infty}
L(x,y_n,p_n) \ge \liminf_{n\to\infty}g_n(p_n) \ge g(p_0)\ge
L(x,y_0,p_0),$ which shows the claim.
\end{proof}

An existence result is an immediate corollary.

\begin{thm}\label{exist-m}
Suppose that $L(x,y,p)$ is a Lagrangian that is convex in $p$ for
each fixed $(x,y)$ and lower semicontinuous in $y$ for each fixed
$(x,p)$. Then for any $a<b$ and $A,B\in\R$ there is a minimizer
$u\in \AC[a,b]$ for which $u(a)=A$ and $u(b)=B$.
\end{thm}

\begin{proof}
Choose $u_n\in\AC[a,b]$ so that
$u_n(a)=A$, $u_n(b)=B$ and
$$\lim_{n\to\infty} {\F}(u_n;a,b) = {\E}(a,A;b,B) < \infty.$$
By Proposition~\ref{equi} the sequence $u_n$
has a subsequence weakly converging to some $u\in\AC[a,b]$
which still satisfies the same boundary conditions,
and by Theorem~\ref{wlsc},
${\F}(u;a,b)\le \lim_{n\to\infty} {\F}(u_n;a,b) = {\E}(a,A;b,B)$,
hence $u$ is a minimizer.
\end{proof}

\subsection{Generalized minimizers}\label{gm}

Here we briefly discuss several possible ways of extending the
notion of minimizers. Initially we consider generalized and
constrained minimizers. Later we give yet another, $\GE$-minimizers,
used for technical purposes only and which, in fact, turn out to be
equivalent to the notion of constrained minimizers. In
Section~\ref{relax}, we describe the relationship between
constrained and relaxed minimizers.

Our original reason for introducing them was to enlarge the
universal singular set to a `generalized universal singular set',
for which our size estimates could still be valid. Although this is
true, we show in Proposition~\ref{universal} that for Lagrangians
satisfying the classical assumptions `generalized universal singular
sets' are in fact universal singular sets. Even for Lagrangians that
are only continuous in $(y,p)$ (and not necessarily convex in $p$),
we will see in Section~\ref{relax} that their `generalized universal
singular sets' are in fact universal singular sets for the
Lagrangian convexified in $p$. Therefore, since it can hardly lead
to any confusion, we skip the `generalized' and, after proving the
necessary results, use the term `universal singular sets' even for
sets defined via generalized minimizers.

We also use generalized minimizers in Section~\ref{tr} to give
non-technical formulations of a variant of Tonelli's regularity theorem,
which shows that Tonelli's results have remarkable stability.

We say that $u\in C[a,b]$
is a \emph{generalized minimizer} for the Lagrangian
\(L\) on \([a,b]\) if its restriction to $(a,b)$
is a locally uniform limit of a sequence
$u_n\in\AC[a_n,b_n]$ such that
$\EE(u_n;a_n,b_n)\to 0$.

We say that $u\in C[a,b]$
is a \emph{constrained minimizer} for the Lagrangian
\(L\) on \([a,b]\) if it
is a uniform limit of a sequence
$u_n\in\AC[a,b]$ such that $u_n(a)=u(a)$, $u_n(b)=u(b)$ and
$\EE(u_n;a,b)\to 0$.

Some remarks may be in order. First, we should explain that
in the definition of generalized minimizers, the
convergence of $u_n$ to $u$ on $(a,b)$ requires that
$\limsup_{n\to\infty} a_n\le a$ and
$b\le\liminf_{n\to\infty} b_n$. However
Lemma~\ref{excess} shows that, equivalently, we could have
required that $a_n\to a$ and $b_n\to b$.
More generally,
if $\eps_n\to 0$
and $a_n,b_n\in (a,b)$ with $a_n\to a$, $b_n\to b$,
then for every generalized minimizer $u$ for \(L\) on \([a,b]\)
there are $u_n\in\AC[a_n,b_n]$ such that
$|u_n-u|<\eps_n$ on $[a_n,b_n]$ and
$\EE(u_n;a_n,b_n)\to 0$: it suffices to replace the
$u_n$ from the definition by $u_{k_n}$ such that
$|u_{k_n}-u|<\eps_n$ on $[a_n,b_n]$ and use Proposition~\ref{excess}.
With such $a_n$, $b_n$ and $u_n$, we use that
$$\F(u_n;a_n,b_n)\le \F(v_n;a_n,b_n)+\EE(u_n;a_n,b_n),$$
where $v_n$ is an affine function joining
$(a_n,u_n(a_n))$ and $(b_n,u_n(b_n))$, together with the fact
that the sequence
$\F(v_n;a_n,b_n)$ is bounded,
to infer from the equiintegrability bound of Proposition~\ref{equi}
that the $u_n$ are equiabsolutely continuous. Hence generalized minimizers
are in fact absolutely continuous. (In what follows, we therefore
consider only functions $u\in\AC[a,b]$.)

The above argument also shows that the same notion of minimizers
would be obtained had we required just pointwise convergence of
$u_n$ to $u$ on $(a,b)$.

It is obvious that constrained minimizers are generalized minimizers.
And
the notion of minimizers has been weakened so much that
an existence result is essentially trivial:
By taking, for any $a<b$ and $A,B\in\R$,
a sequence $u_n\in\AC[a,b]$ such that $u_n(a)=A$,
$u_n(b)=B$ and $\EE(u_n;a,b) \to 0$, and using Proposition~\ref{equi}
together with the Arzela-Ascoli Theorem,
we see that constrained minimizers exist. We record some of these
facts in the following theorem.

\begin{thm}\label{exi}
Let $L$ be a Lagrangian.
Suppose that $a<b$ and $A,B\in\R$. Then there is a constrained minimizer
such that $u(a)=A$ and $u(b)=B$. Moreover, every minimizer is a constrained
minimizer, every constrained minimizer
is a generalized minimizer and every generalized minimizer is
absolutely continuous.
\end{thm}

Notice also that Lemma~\ref{excess} implies that the restriction
of a generalized minimizer to a subinterval is a generalized minimizer;
the same statement for constrained minimizers is not so obvious (but
true; see Corollary~\ref{constrained-prop}). As with ordinary
minimizers, the maximum and minimum of two constrained minimizers with the
same boundary conditions is a constrained minimizer.

\begin{prp}\label{gm-min}
Let $u,v$ be constrained minimizers for $L$ on $[a,b]$ such that
$u(a)\le v(a)$ and $u(b)\le v(b)$. Then $\max(u,v)$
and $\min(u,v)$ are
constrained minimizers for $L$ on $[a,b]$.
\end{prp}

\begin{proof}
Let $u_n,v_n\in\AC[a,b]$ be such that $u_n(a)=u(a)$, $v_n(a)=v(a)$,
$u_n(b)=u(b)$, $v_n(b)=v(b)$,
$\EE(u_n;a,b) \to 0$ and $\EE(v_n;a,b) \to 0$.
Then
$\EE(\max(u_n,v_n);a,b)+\EE(\min(u_n,v_n);a,b)
=\EE(u_n;a,b)+\EE(u_n;a,b)\to 0,$
hence $\EE(\max(u_n,v_n);a,b)\to 0$ and $\EE(\min(u_n,v_n);a,b)\to 0$.
\end{proof}

For generalized minimizers this argument fails, and indeed,
the maximum of two generalized minimizers with the
same boundary conditions need not be a generalized minimizer.
A variant of the argument shows that
the analogue of the
above statement holds for generalized minimizers $u$ and $v$
provided we assume that
$u(a)< v(a)$ and $u(b) < v(b)$. (These remarks are not used
in this paper.)

Easy examples show that none of the inclusions from
the last sentence of Theorem~\ref{exi} may be reversed.
Generalized minimizers may fail to be constrained even
for Lagrangians satisfying the classical assumptions
(under which, of course, constrained minimizers coincide with minimizers);
for example for
$$L(x,y,p)= (x^3 - y^5)^2 p^{20} + \eps_0 p^2,$$
where $\eps_0>0$ is a small enough constant. Indeed, let
$u_n\in\AC[0,1]$ be a minimizer for $u_n(0)=-1/n$, $u_n(1)=1$. The
usual estimates proving the Lavrentiev phenomenon (see
either~\cite{lavr-phen} or~\cite{butt_gia,Lavrentiev,mania}) show
that there is $\delta>0$ such that $u_n(x)\le \tfrac{1}{2}x^{3/5}$
for all $n$ and all $0\le x\le\delta$. Hence the limit $u$ of (a
subsequence of) the $u_n$'s stays below $\tfrac{1}{2} x^{3/5}$ on
$[0,\delta]$. But, by the Lavrentiev phenomenon estimates mentioned
above, any such $u$ is of energy larger than that of $\tilde
u(x)=x^{3/5}$ and so is only a generalized minimizer, not a
minimizer.

In this example, we see that a generalized minimizer fails to
be constrained because it is on a higher energy level. We now show that
this is the only way this phenomenon can occur. This will then lead us
to recognition that our two notions of generalized minimizers give rise to
the same universal singular set.
For this, it is convenient to introduce some notation
for the `generalized energy'.
For $u\in\AC[a,b]$ denote by $\GF(u;a,b)$ the infimum of all possible
$\liminf_{n\to\infty}\F(u_n;a_n,b_n)$ where $a_n\to a$,
$b_n\to b$ and $u_n\in\AC[a_n,b_n]$ converge to $u$ on $(a,b)$. We also let
$$\GE(a,A;b,B)=
\inf\left\{\GF(u;a,b): u\in \AC [a,b],\, u(a)=A,\, u(b)=B\right\}.
$$
and say that $u$ is an \emph{$\GE$-minimizer} (for $L$ on $[a,b]$) if
$$\GF(u;a,b)=\GE(a,u(a);b,u(b)).$$

\begin{lem}\label{why not}
For any $a<b$ and $A,B\in\R$,
$\GE$-minimizers exist.
\end{lem}

\begin{proof}
By Proposition~\ref{equi}
and the Arzela-Ascoli Theorem, it suffices to find
$u_n\in\AC[a_n,b_n]$ such that $(a_n,u_n(a_n))\to(a,A)$,
$(b_n,u_n(b_n))\to(b,B)$ and
$$\limsup_{n\to\infty}\F(u_n;a_n,b_n)\le \GE(a,A;b,B).$$
To find the $u_n$ we start by choosing $v_n\in \AC [a,b]$ such that
$$v_n(a)=A,\; v_n(b)=B\textrm{ and }\GF(v_n;a,b) < \GE(a,A;b,B)+1/n.$$
We then choose $a<a_n<a+1/n$ and $b-1/n<b_n<b$ such that
$|v_n(a_n)-v_n(a)|<1/n$ and $|v_n(b_n)-v_n(b)|<1/n$
and finish by using the definition of $\GF(v_n;a,b)$ to find
$\alpha_n\in (a -1/n,a_n)$, $\beta_n\in(b_n, b+1/n)$,
and $u_n\in\AC[\alpha_n,\beta_n]$ such that
$|u_n(a_n)-v_n(a_n)|<1/n$, $|u_n(b_n)-v_n(b_n)|<1/n$ and
$\F(u_n;\alpha_n,\beta_n)<\GF(u_n;a,b)+1/n.$
Then
$\limsup_{n\to\infty}\F(u_n;a_n,b_n)\le
\limsup_{n\to\infty}\F(u_n;\alpha_n,\beta_n)\le
\lim_{n\to\infty}\GF(u_n;a,b)=\GE(a,A;b,B).$
\end{proof}

\begin{lem}\label{lem-constrained}
Every $\GE$-minimizer $u\in\AC[a,b]$
is a (uniform) limit of a sequence $u_n\in \AC[a,b]$ such that
$u_n(a)=u(a)$, $u_n(b)=u(b)$ and $\F(u_n;a,b)\to \GF(u;a,b).$
\end{lem}

\begin{proof}
Let $\eps>0$; we have to find a $v\in\AC[a,b]$ such that
$v(a)=u(a)$, $v(b)=u(b)$,
$|v(x)-u(x)|<\eps$ on $[a,b]$ and
$\F(v;a,b)< \GF(u;a,b)+\eps.$

Find $v_n\in\AC[a_n,b_n]$ converging to $v$ on $(a,b)$ such that
$\F(v_n;a_n,b_n)< \GF(u;a,b) + 2^{-n-4}\eps.$
For $a<s<t<b$ denote
$\phi(s,t)=\limsup_{n\to\infty} \F(v_n;s,t).$
We shall assume that for rational $s,t$ this $\limsup$
is actually a limit; if necessary, this can be achieved by
passing to a subsequence of $v_n$.

We claim that for every $\eta>0$, there is $\delta>0$ such that
$\phi(s,t)<\eta$, whenever $a<s<t<a+\delta$.
To see this, note that $\sigma:=\sup_{a<s<t<b} \phi(s,t)
\le \GF(u;a,b)<\infty$ and find $a<s_0<t_0<b$ such that
$\phi(s_0,t_0)>\sigma-\eta$. If
$\phi(s,t)\ge \eta$ for some $a<s<t<s_0$,
we would pick rational $\tilde s\in(a,s)$ and
$\tilde t\in(t,t_0)$, and use that
$\phi(\tilde s,\tilde t)=\lim_{n\to\infty} \F(v_n;\tilde s,\tilde t)$
to obtain a contradiction by estimating
$\phi(\tilde s,t_0)\ge\phi(\tilde s,\tilde t)+\phi(s_0,t_0) >\sigma$.

Similarly we show that for every $\eta>0$ there is $\delta>0$ such that
$\phi(s,t)<\eta$ whenever $b-\delta<s<t<b$.

Hence, for each $j=1,2,\dots$ we may find $\delta_j>0$ such that
$\phi(s,t)<2^{-j-4}\eps$ provided that either $a<s<t<a+\delta_j$
or $b-\delta_j<s<t<b$.

Choose $\alpha_j\searrow a$ and $\beta_j\nearrow b$
so that $\alpha_j<a+\delta_j$, $\beta_j>b-\delta_j$,
$\alpha_1<\beta_1$ and $u$ has a finite derivative at each of
the points $\alpha_j,\beta_j$. Choose also
$0<\tau_j<\min(\alpha_j-\alpha_{j+1},\beta_{j+1}-\beta_j)$
and $C_j<\infty$
so that
$|u(x)-u(\alpha_j)|\le C_j|x-\alpha_j|$ whenever
$|x-\alpha_j|\le\tau_j$
and
$|u(x)-u(\beta_j)|\le C_j|x-\beta_j|$ whenever
$|x-\beta_j|\le\tau_j$.

Let $R_j:=(C_j+2)(b-a)+\max_{t\in[a,b]} |u(t)|$ and
find $M_j<\infty$ such that $L(x,y,p)\le M_j$ for
$x\in[a,b]$, $|y|\le R_j$ and
$|p|\le C_j+2$. Also find $0<\eta_j<\tau_j$ so that
$\eta_j M_j<2^{-j-2}$.

For each $j$ choose $n_j\ge j$ so large that
$a_{n_j}<\alpha_j$, $b_{n_j}>\beta_j$,
$|u_{n_j}(x)-u(x)|<2^{-j-4}\eta_j$ on $[\alpha_{j},\beta_j]$,
and, if $j\ge 2$,
$\F(v_{n_{j}}; \alpha_{j},\alpha_{j-1})<2^{-j-3}\eps$
and $\F(v_{n_{j}}; \beta_{j-1},\beta_{j})<2^{-j-3}\eps$.
(To get the last requirements we have used that
$\phi(\alpha_{j},\alpha_{j-1})<2^{-j-3}\eps$
and $\phi(\beta_{j-1},\beta_{j})<2^{-j-3}\eps$.)

Define $v\colon [a,b]\to\R$ by
$$
v(x)=
\begin{cases}
u(x) & \textrm{for $x=a$ and $x=b$}\\
v_{n_1}& \textrm{on  $[\alpha_1,\beta_1]$}\\
v_{n_{j+1}}(x) & \textrm{on  $[\alpha_{j+1},\alpha_{j}-\eta_j]
\cup [\beta_{j}+\eta_j,\beta_{j+1}]$}\\
\textrm{affine}& \textrm{on $[\alpha_{j}-\eta_j,\alpha_{j}]$
and on $[\beta_{j}, \beta_{j}+\eta_j]$.}
\end{cases}
$$

Clearly, $v$ is locally absolutely continuous on $(a,b)$
and, since $|v(x)-u(x)|<\eta_j$ for
$x\in(a,\alpha_j)\cup(\beta_j,b)$, it is
also continuous on $[a,b]$.

Since on $(\alpha_{j}-\eta_j,\alpha_{j})$,
$|v|\le\max(v_{n_j}(\alpha_j),v_{n_{j+1}}(\alpha_{j}-\eta_j)|)+(b-a)\le
R_j$ and
$|v'|=|v_{n_j}(\alpha_j)-v_{n_{j+1}}(\alpha_{j}-\eta_j)|/\eta_j \le
|u(\alpha_j)-u(\alpha_{j}-\eta_j)|/\eta_j + 2\le C_j+2$, we have
$\F(v; \alpha_{j}-\eta_j,\alpha_{j})\le M_j\eta_j<2^{-j-4}\eps$. A
similar argument gives that $\F(v; \beta_{j}, \beta_{j}+\eta_j)\le
M_j\eta_j<2^{-j-4}\eps$. Also, recall that $n_{j+1}$ has been chosen
so that $\F(v_{n_{j+1}};
\alpha_{j+1},\alpha_{j}-\eta_j)<2^{-j-4}\eps$ and $\F(v_{n_{j+1}};
\beta_{j}+\eta_j,\beta_{j+1})<2^{-j-4}\eps$, and infer from $L\ge 0$
that $\F(v_{j_1},\alpha_1,\beta_1)\le \F(v_{j_1},a_1,b_1)<
\GF(u;a,b) + \eps/2$.

Adding all these estimates of $\F$ together gives
$\int L(x,v(x),v'(x))\,dx <\GF(u;a,b) + \eps.$
Since this also implies that $v'$ is integrable over $[a,b]$,
we see that $v$, being an indefinite Lebesgue integral of
$v'$ on $[a,b]$, is absolutely continuous on $[a,b]$.
Hence $v$ has all the required properties and we are done.
\end{proof}

\begin{cor}\label{constrained}
$\GE(a,A;b,B)=\E(a,A;b,B)$
for every $a,A,b,B$, and
the notions of constrained and $\GE$-minimizers coincide.
\end{cor}

\begin{proof}
Clearly, $\GE(a,A;b,B)\le \E(a,A;b,B)$.
Hence, if $u$ is any $\GE$-minimizer, Lemma~\ref{lem-constrained}
shows that it is also a constrained minimizer
and that $\E(a,A;b,B)\le \GE(a,A;b,B)$; since $\GF$ minimizers exist,
this also shows that $\GE(a,A;b,B)=\E(a,A;b,B)$. Having this equality,
it is clear that every constrained minimizer is an $\GE$-minimizer.
\end{proof}

Various properties of constrained minimizers follow from these facts.
We just record the following

\begin{cor}\label{constrained-prop}
The restriction of a constrained minimizer to a
subinterval is a constrained minimizer.
\end{cor}

\begin{proof}
Let $u$ be a constrained minimizer on $[a,b]$,
with the corresponding $u_n\in\AC[a,b]$ and let $c\in(a,b)$.
Denoting by $v_n$ and $w_n$ the restrictions of $u_n$
to $[a,c]$ and $[c,b]$, respectively,
and using Corollary~\ref{constrained},
we find
\begin{align*}
\E(a,u(a);b,u(b)) &\le
\E(a,u(a);c,u(c))+\E(c,u(c);b,u(b))\\
&=\GE(a,u(a);c,u(c))+\GE(c,u(c);b,u(b))\\
&\le \GF(u;a,c)+\GF(u;c,b)\\
&\le \liminf_{n\to\infty}
(\F(v_n;a,c)+\F(w_n;c,b))\\
&=\E(a,u(a);b,u(b)).
\end{align*}
Since $\GF(u;a,c)\ge\GE(a,u(a);c,u(c))$ and
$\GF(u;c,b)\ge \GE(c,u(c);b,u(b))$,
we have
$$\GF(u;a,c)=\GE(a,u(a);c,u(c))\textrm{ and }\GF(u;c,b)= \GE(c,u(c);b,u(b)).$$
Hence the restrictions of $u$ to $[a,c]$ and $[c,b]$ are
$\GE$-minimizers, and thus
constrained minimizers by Corollary~\ref{constrained}.
\end{proof}

\begin{lem}\label{lem-universal}
Suppose that $u\in\AC[a,b]$ is a generalized minimizer for $L$
such that $u'_+(a)=\infty$. Then for every $B>u(b)$ there is
a constrained minimizer $v$ on $[a,b]$ such that $v(a)=u(a)$, $v(b)=B$, and
$v$ has infinite right derivative at $a$.
\end{lem}

\begin{proof}
Let $w\in\AC[a,b]$ be a constrained minimizer
such that $w(a)=u(a)$ and $w(b)=B$. If the right derivative of
$w$ at $a$ is infinite, we are done. So assume that this is not the case,
hence $\liminf_{x\searrow a}\frac{|w(x)-w(a)|}{|x-a|} < \infty.$
Since $u'_+(a)=+\infty$, this implies, in particular, that there is
$c\in(a,b)$ such that $w(c)<u(c)$. Also use that $w(b)>u(b)$
to find $c<d<b$ so that $w(d)>u(d)$.

Let $c>\beta_n\searrow a$ and $M<\infty$ be such that
$|w(\beta_n)-w(a)|< M (\beta_n-a)$ for each $n$,
and find $a<\alpha_n<\beta_n$ so that
$|w(\beta_n)-w(\alpha_n)|< M (\beta_n-\alpha_n).$
Denoting $\delta_n=\tfrac12\min(\beta_n-\alpha_n,b-d,u(c)-w(c),B-u(b)),$
we use that $w$ is a constrained minimizer to find
$w_n\in\AC[a,b]$ such that
$w_n(a)=w(a)$, $w_n(b)=w(b)$,
$|w_n(x)-w(x)|<\delta_n$ on $[a,b]$, and
$\F(w_n;a,b)< \GF(w;a,b) + 1/n.$
Also, we use that $u$ is a generalized minimizer to find
$a_n<\alpha_n<b-\delta_n <b_n$ and
$u_n\in\AC[a_n,b_n]$ such that
$|u_n(x)-u(x)|<\delta_n$ on $[\alpha_n,b-\delta_n]$ and
$\EE(u_n;a_n,b_n)< 1/n.$

Since $w_n(c)<u_n(c)$ and $w_n(d)>u_n(d)$,
there are $c_n\in(c,d)$ such that
$u_n(c_n)=w_n(c_n)$.
By considering the function defined as $u_n$ on
$[a_n,a+\delta_n]\cup[c_n,b_n]$, as $w_n$ on $[\alpha_n,c_n]$
and affine on $[a+\delta_n,\alpha_n]$, we see that
$\F(u_n;a_n,b_n)\le
\F(u_n;a_n,a+\delta_n)+\F(u_n;c_n,b_n)+\F(w_n;\alpha_n,c_n)
+M(\beta_n-\alpha_n)+\EE(u_n;a_n,b_n),$
hence
$\F(u_n;\alpha_n,c_n)\le\F(w_n;\alpha_n,c_n)+2/n.$

Defining $v_n$ as $u_n$ on $[\alpha_n,c_n]$
and as $w_n$ on $[c_n,b_n]$, we therefore have
$\F(v_n;\alpha_n,b)
=\F(u_n;\alpha_n,c_n)+\F(w_n;c_n,b)
\le \F(w_n;\alpha_n,b)+2/n
\le \F(w_n;a,b)+2/n
\to \GF(w;a,b).$
Hence, choosing a subsequence along which $c_n$ converges to,
say, $\tilde c\in[c,d]$, and defining $v$ as $u$ on $[a,\tilde c]$
and as $w$ on $[\tilde c,b]$, we see that
$\GF(v;a,b)\le \GF(w;a,b)=\GE(a,v(a);b,v(b)).$
So $v$ is an $\GE$-minimizer, and so, by Corollary~\ref{constrained},
it is a constrained minimizer. Finally,
$v'_+(a)=\infty$ since $v\ge u$ on $[a,c]$.
\end{proof}

\begin{prp}\label{universal}
The following statements about a point $(x_0,y_0)\in\R^2$
are equivalent.
\begin{enumerate}
\item\label{universal.1}  There
are $a\leq  x_0\leq b$ with \(a<b\) and a generalized
minimizer $u$ for \(L\) on \([a,b]\) such that
$u(x_0)=y_0$ and $u$ has infinite derivative at $x_0$.
\item\label{universal.2} There
are $a\leq  x_0\leq b$ with \(a<b\) and a constrained
minimizer $u$ for \(L\) on \([a,b]\) such that
$u(x_0)=y_0$ and $u$ has infinite derivative at $x_0$.
\end{enumerate}
Moreover, if $L(x,y,p)$ is continuous in $(y,p)$ and
convex in $p$, these statements hold if and only if
\begin{enumerate}\setcounter{enumi}{2}
\item\label{universal.3}  There
are $a\leq  x_0\leq b$ with \(a<b\) and a
minimizer $u$ for \(L\) on \([a,b]\) such that
$u(x_0)=y_0$ and $u$ has infinite derivative at $x_0$.
\end{enumerate}
\end{prp}

\begin{proof}
By Lemma~\iref{lem-universal}, \iref{universal.2} holds
provided that there is $b>x_0$ and a generalized
minimizer $u$ for \(L\) on \([a,b]\) such that
$u(x_0)=y_0$ and $u'_+(x_0)=\infty$.
Symmetric arguments shows that \iref{universal.2} also holds
if $u'_-(x_0)=-\infty$ or if
there is $a<x_0$ and a generalized
minimizer $u$ for \(L\) on \([a,b]\) such that
$u(x_0)=y_0$ and $u$ has infinite left derivative at $x_0$.
Hence \iref{universal.1}$\Rightarrow$\iref{universal.2}.

The implication \iref{universal.2}$\Rightarrow$\iref{universal.1}
is obvious since every constrained minimizer is a generalized one.

Under the additional assumptions on the Lagrangian, \iref{universal.2}
is equivalent to \iref{universal.3}, since then the notions of
constrained minimizers and minimizers coincide.
\end{proof}

This statement allows us to define
the \emph{universal singular set of a Lagrangian} as
all points \( (x_0,y_0)\in\R^2\) for which the statement
\iref{universal.1}, or equivalently \iref{universal.2}, of
Proposition~\ref{universal} holds.

We now return to our example illustrating the difference between
general and constrained minimizers and observe that not only were
the two minimizers on different energy levels, but  behind the whole
discrepancy was the fact that an end-point of the minimizers
belonged to the universal singular set:

\begin{prp}\label{bk2}
Suppose that $u\in\AC[a,b]$ is
such that
$$ \liminf_{x\searrow a}\frac{|u(x)-u(a)|}{|x-a|} < \infty
\textrm{ and }
\liminf_{x\nearrow b}\frac{|u(x)-u(b)|}{|x-b|} < \infty
$$
and that $v_n$ is a sequence of absolutely continuous functions
on $[a_n,b_n]$ such that
$(a_n,v_n(a_n))\to (a,u(a))$, $(b_n,v_n(b_n))\to (b,u(b))$
and $\EE(v_n;a_n,b_n)\to 0$ (but possibly converging to
a limit different from $u$).
Then $$\limsup_{n\to\infty} \F(v_n;a_n,b_n)\le\GF(u;a,b).$$
\end{prp}

\begin{proof}
We may assume that $\GF(u;a,b)<\infty$ and
$\lim_{n\to\infty} \F(v_n;a_n,b_n)$ exists;
this will allow us to pass to a subsequence of $v_n$.
Denote by $M$ the maximum of the two
limits inferior in the assumption and choose
$\alpha_n\searrow a$ and $\beta_n\nearrow b$
such that $\alpha_n<\beta_n$,
$|u(\alpha_n)-u(a)|\le (M+1)|\alpha_n-a|$
and $|u(\beta_n)-u(a)|\le (M+1)|\beta_n-b|$.
Let $\delta_n=\min(\alpha_n-a,b-\beta_n)/2$.
Passing to a subsequence of $v_n$ if necessary,
we also assume
$\|(a_n,v_n(a_n)) - (a,u(a))\|<\delta_n$ and
$\|(b_n,v_n(b_n))- (b,u(b))\|<\delta_n$.

Find $u_n\in\AC[\alpha_n,\beta_n]$ such that
$|u_n(x)-u(x)|<\delta_n$ on $[\alpha_n,\beta_n]$
and $\GF(u;a,b)=\lim_{n\to\infty} \F(u_n;\alpha_n,\beta_n)$.

Since $a_n <\alpha_n<\beta_n < b_n$, we may define
$w_n$ on $[a_n,b_n]$ which agrees with $v_n$ at the points $a_n,b_n$,
with $u_n$ on $[\alpha_n,\beta_n]$ and is affine
on $[a_n,\alpha_n]$ and $[\beta_n,b_n]$.

Choose $R>0$ so large that
the graphs of $v_n$, $u_n$, and so also of $w_n$ are contained in
$[-R,R]^2$. Let $C$ be an upper bound for $L(x,y,p)$
for $|x|,|y|\le R$ and $|p|\le M+4$.

Observing that $\alpha_n-a_n\ge\delta_n$, we estimate that
$|w_n(a_n)-w_n(\alpha_n)|\le |u(a)-u(\alpha_n)|+2\delta_n
\le (M+4)|a-\alpha_n|
$
and, similarly, that
$|w_n(b_n)-w_n(\beta_n)|\le (M+4)|b_n-\beta_n|$.
Hence, letting $l_n=|\alpha_n-a|+|b-\beta_n|$,
\begin{align*}
\F(v_n;a_n,b_n)
&\le \F(w_n;a_n,b_n)+\EE(v_n;a_n,b_n)\\
&\le \F(w_n;\alpha_n,\beta_n) +Cl_n+\EE(u_n;a_n,b_n)\\
&\to \GF(u;a,b),
\end{align*}
and it suffices to take a limit.
\end{proof}

Using the definition of $\GF(v;a,b)$ and
Proposition~\ref{constrained}, we have the following result.

\begin{cor}\label{bk3}
If $u$ satisfies the assumptions of Proposition~\ref{bk2}
and $v$ has the same boundary values then
$\GF(v;a,b)\le\GF(u;a,b).$

If $u$ is a constrained minimizer
satisfying the assumptions of Proposition~\ref{bk2},
then every generalized minimizer having the same boundary values
is constrained.
\end{cor}

\subsection{Approximation and Relaxation}\label{relax}

In this section, we show how the usual relaxation (i.e.\
convexification) procedure applies to the study of universal
singular sets for general Lagrangians. A special case of
Theorem~\ref{generelaxed} says that, under the additional assumption
that the Lagrangian $L$ is also continuous in $(y,p)$, generalized
minimizers of $L$ starting and ending at points outside the
universal singular set of the convexified Lagrangian $L^c$ are
necessarily minimizers of $L^c$. As a corollary, in
Theorem~\ref{univrelaxed} we show that under the same assumptions on
$L$, the generalized universal singular set of $L$ coincides with
the universal singular set of its convexification \(L^c\).

All Lagrangians $L(x,y,p)$
in this section are assumed to be continuous in
$(y,p)$, for each fixed $x$.

Before introducing relaxation, we show that all $u\in AC[a,b]$ can
be approximated (in the space $W^{1,1}(a,b)$, and whilst preserving
the boundary conditions) by a function whose energy is not much
greater than that of $u$ and which is $C^1$ on a dense open subset
of $(a,b)$ of full measure. This would be easy if $u$ were Lipschitz
on a dense open subset of $(a,b)$ of full measure; however, under
our assumptions, this need not be the case.

\begin{prp}\label{approx1}
For every $u\in\AC[a,b]$ with $\F(u;a,b)<\infty$ and every $\eps>0$,
there is $v\in\AC[a,b]$ such that
$$v(a)=u(a), \, v(b)=u(b),\,
\int_a^b|v'-u'|\,dx<\eps,\, \F(v;a,b)\le \F(u;a,b)+\eps,$$
and
almost every point of $[a,b]$ has a neighbourhood on which $v$ is
affine.
\end{prp}

\begin{proof}
For $a\le\alpha<\beta\le b$ let $u_{\alpha,\beta}$ denote the affine
function that joins $(\alpha,u(\alpha))$ and $(\beta,u(\beta))$. Let
$\eta=\eps/(b-a)$.

Our proof depends upon showing that for almost every $\alpha\in
(a,b)$,
\begin{equation}\label{approx1.1} \F(u_{\alpha,\beta};\alpha,\beta)
< \F(u;\alpha,\beta) +\eta (\beta-\alpha),
\end{equation} whenever
$\beta>\alpha$ is close enough to $\alpha$.

We start by showing how~(\ref{approx1.1}) gives the Proposition: Let
$\mathcal{J}$ be the family of intervals $J=[\alpha,\beta]\subset
(a,b)$ such that
$$
\int_\alpha^\beta|u_{\alpha,\beta}'(x)-u'(x)|\,dx<\eta |J|
\textrm{ and }
\F(u_{\alpha,\beta};\alpha,\beta)
<
\F(u;\alpha,\beta)+\eta |J|
$$
Since almost every $\alpha$ is a Lebesgue point of $u'$ for
which~(\ref{approx1.1}) holds, almost every $\alpha$ has the
property that $[\alpha,\beta]\in\mathcal{J}$ whenever $\beta>\alpha$
is close enough to $\alpha$. Hence $\mathcal{J}$ covers $(a,b)$ in
the sense of Vitali, and so the Vitali covering theorem provides us
with disjoint $J_j=[\alpha_j,\beta_j]\in\mathcal{J}$ covering almost
all of $(a,b)$.

Define $v\colon [a,b]\to\R$ by
  \[v(x)=\begin{cases}
    u(x), &\textrm{for }x\in[a,b]\setminus\bigcup_j J_j,\\
    u_{\alpha_j,\beta_j}(x), &\textrm{on }J_j.
  \end{cases}\]
Since $\int_{J_j}|v'-u'|\,dx<\eta |J_j|$ and $v$ agrees with $u$ at
the endpoints of the intervals $J_j$ and outside their union, we see
that $v'$ is integrable, $\int_a^b|v'-u'|\,dx<\eta(b-a)=\eps$ and
$v$ is an indefinite integral of $v'$. Consequently, $v\in\AC[a,b]$
and, since $\{J_j\}$ covers almost all of $(a,b)$, $v$ is affine in
a neighbourhood of almost every point and
$$\F(v;a,b) \le\sum_j \F(v;\alpha_j,\beta_j)
\le \sum_j (\F(u;\alpha_j,\beta_j)+\eta|J_j|)
<
\F(u;a,b) +\eps.$$

It remains to prove that~(\ref{approx1.1}) holds for almost every
\(\alpha\in (a,b)\): For those $x\in(a,b)$ at which $u$ is
differentiable, define $g(x)$ to be the largest number for which
$$|L(x,y,p)-L(x,u(x),u'(x))|<\tfrac{1}{2}\eta\textrm{, whenever }
                \|(y,p)-(u(x),u'(x))\|<g(x).$$
Since $L(x,y,p)$ is a continuous function in $(y,p)$ for every
\(x\), $g$ is a strictly positive measurable function.
(Measurability follows upon observing that $\{x: g(x)<c\}$ is the
union of the sets $\{x:
|L(x,y,p)-L(x,u(x),u'(x))|>\tfrac{1}{2}\eta\}$ over rational $y,p$
for which $\|(y,p)-(u(x),u'(x))\|<c$.)

Let $\alpha\in(a,b)$ be a point at which $u$ is differentiable and
at which both $g$ and $u'$ are approximately continuous. Since both
the values and the slopes
 of $u_{\alpha,\beta}$ have a
bound independent of $\beta\in(\alpha,b)$, there is $C<\infty$ such
that $L(x,u_{\alpha,\beta}(x),u_{\alpha,\beta}'(x))\le C$ for every
$\beta>\alpha$ and $x\in(\alpha,\beta)$. If $\beta>\alpha$ is close
enough to $\alpha$, then
$|u_{\alpha,\beta}(x)-u(x)|<\tfrac{1}{2}g(\alpha)$ for every
$x\in(\alpha,\beta)$ and the set $$T=\{x\in(\alpha,\beta):
|u_{\alpha,\beta}'(x)-u'(x)|>\tfrac{1}{2}g(\alpha) \textrm{ or }
|g(x)-g(\alpha)|>\tfrac{1}{2}g(\alpha)\}$$ has measure less that
$\eta(\beta-\alpha)/(2C)$. Letting $S=(\alpha,\beta)\setminus T$ and
observing that for almost all $x\in S$,
$\|(u_{\alpha,\beta}(x),u_{\alpha,\beta}'(x))-(u(x),u'(x))\|<g(\alpha)$
and so $|L(x,u_{\alpha,\beta}(x),u_{\alpha,\beta}'(x)) -
L(x,u(x),u'(x))|<\tfrac{1}{2}\eta$, we find
\begin{align*}
\F(u_{\alpha,\beta};\alpha,\beta) &=\int_S
L(x,u_{\alpha,\beta}(x),u_{\alpha,\beta}'(x))\,dx
+\int_T L(x,u_{\alpha,\beta}(x),u_{\alpha,\beta}'(x))\,dx\\
&\le
\int_S L(x,u(x),u'(x))\,dx + \tfrac{1}{2}\eta|S| + C|T|\\
&< \F(u;\alpha,\beta) +\eta (\beta-\alpha),
\end{align*}
which finishes the proof.
\end{proof}

We now come to the main part of this section. Since we will work
with two Lagrangians simultaneously, we write $\F_L(u;a,b)$ instead
of $\F(u;a,b)$ for $\int_a^b L(x,u(x),u'(x)\,dx$;  we use
$\E_L(a,A;b,B)$ and $\EE_L(u;a,b)$ similarly. We define relaxed
minimizers as follows. Given a Lagrangian $L$, we denote by $L^c$
the convexification of $L$ with respect to the third variable; thus
$L^c(x,y,p)$ is equal to
  $$\inf\{ \lambda_1 L(x,y,p_1)+ \lambda_2L(x,y,p_2) :
\lambda_i\ge 0,\,\lambda_1+ \lambda_2=1 \textrm{ and } p=\lambda_1
p_1 +\lambda_2 p_2 \}.$$
 We say that $u\in AC[a,b]$ is a
\emph{relaxed minimizer for $L$} if it minimizes the $L^c$ energy
with respect to its boundary data in $a$ and $b$, that is,
$$\F_{L^c}(u;a,b)=\E_{L^c}(a,u(a);b,u(b)).$$

Provided that the (superlinear) lower bound $\omega$ for $L$ is
convex, it is clear that $L^c$ has the same lower bound.

\begin{lem}
  If \(L\) be a Lagrangian that is continuous in \( (y,p)\) for each fixed \(x\), then
    $L^c$ is continuous in $(y,p)$ for each fixed \(x\).
\end{lem}

\begin{proof}
Let $x$ be fixed, and suppose $(y^n,p^n)\to(y^0,p^0)$ and $\eps>0$.
If $\lambda_i\ge 0$ are such that
  $$\lambda_1+ \lambda_2=1,\; p^0=\lambda_1 p_1 +\lambda_2
  p_2,$$
and
  $$\lambda_1 L(x,y^0,p_1)+ \lambda_2 L(x,y^0,p_2) <
  L^c(x,y^0,p^0)+\eps,$$
then we may change the $p_i$ and $\lambda_i$ slightly, relabelling
them if necessary, to get $p_1<p^0<p_2$. Then for large $n$, there
are $\lambda_i^n\ge 0$, with $\lambda_1^n+ \lambda_2^n=1$ for which
$p^n=\lambda_1^n p_1 +\lambda_2^n p_2$; moreover
$\lambda_i^n\to\lambda_i$. Hence
\begin{align*}
  L^c(x,y^n,p^n) & \le \lambda_1^n L(x,y^n,p_1)+ \lambda_2^n
L(x,y^n,p_2)\\
  & \to \lambda_1 L(x,y^0,p_1)+ \lambda_2 L(x,y^0,p_2)\\
  &<  L^c(x,y^0,p^0)+\eps,
\end{align*}
and we see that $\limsup_{n\to\infty}
L^c(x,y^n,p^n)\le L^c(x,y^0,p^0).$

For the opposite direction, it suffices to show that there are
bounded $p_i^n$ and $\lambda_i^n\ge 0$ for which
  $$\lambda_1^n+
  \lambda_2^n=1,\; p^n=\lambda_1^n p_1^n +\lambda_2^n p_2^n$$
 and
  $$\lambda_1^n L(x,y^n,p_1^n)+ \lambda_2^n L(x,y^n,p_2^n) \le
  L^c(x,y^n,p^n)+\eps;$$
the inequality $\liminf_{n\to\infty} L^c(x,y^n,p^n)\ge
L^c(x,y^0,p^0)$ then follows by taking limits over subsequences and
arbitrary $\eps>0$. If $L^c(x,y^0,p^0)>L^c(x,y,p)$, then the
superlinearity of $L$ implies that any $p_i^n$ satisfying the above
conditions are bounded. If $L^c(x,y^0,p^0)=L^c(x,y,p)$, we take
$p_1^n=p_1<p^0<p_2=p_2^n$, sufficiently close to $p^0$; the desired
inequality follows by continuity.
\end{proof}

\begin{lem}\label{relax1}
For every $u\in\AC[a,b]$ and $\eps>0$, there is $v\in\AC[a,b]$ such
that
  $$v(a)=u(a),\; v(b)=u(b),\; |v(x)-u(x)|< \eps \textrm{ for every }
  x\in[a,b],$$
 and
$$\F_L(v;a,b)\le \F_{L^c}(u;a,b)+\eps.$$
\end{lem}

\begin{proof}
Since the statement is obvious when $\F_{L^c}(u;a,b)=\infty$ and
since adding a constant to $L$ does not change the inequalities to
be proved, we  assume that $\F_{L^c}(u;a,b)<\infty$ and
$L(x,y,p)\ge|p|$. By Lemma~\ref{approx1} used for the Lagrangian
$L^c$ we may also assume that almost every point of $[a,b]$ has a
neighbourhood on which $u$ is affine.

Our proof depends upon showing that, given any $\eta>0$, almost
every $\alpha\in (a,b)$ has the property that  for every $\beta\in
(\alpha,\beta )$ that is close enough to $\alpha$, there are
functions $u_{\alpha,\beta}\in\AC[\alpha,\beta]$, for which
$u_{\alpha,\beta}(\alpha)=u(\alpha)$,
$u_{\alpha,\beta}(\beta)=u(\beta)$ and
\begin{equation}
\F_L(u_{\alpha,\beta};\alpha,\beta) < \F_{L^c}(u;\alpha,\beta) +\eta
(\beta-\alpha).\label{relax2}
\end{equation}

We first show how this gives the Lemma. Let $\eta=\eps/(b-a)$ and
let $\mathcal{J}$ be the family of those intervals
$J=[\alpha,\beta]\subset (a,b)$ for which $\eta |J|<\eps/3$,
$\F_{L^c}(u;\alpha,\beta))<\eps/6$, and there is
$v_J\in\AC[\alpha,\beta]$ with $v_J(\alpha)=u(\alpha)$,
$v_J(\beta)=u(\beta)$ and $ \F_L(v_J;\alpha,\beta) <
\F_{L^c}(u;\alpha,\beta) +\eta|J|. $

Since~(\ref{relax2}) implies that $\mathcal{J}$ covers $(a,b)$ in
the sense of Vitali, we may use the Vitali covering theorem to find
disjoint $J_j=[\alpha_j,\beta_j]\in\mathcal{J}$ whose union covers
almost all of $(a,b)$. Let $v_j=v_{J_j}$.

Define $v\colon [a,b]\to\R$ by
  \[v(x)=\begin{cases}
            u(x), &\textrm{for }x\in[a,b]\setminus\bigcup_j J_j,\\
            v_j(x), &\textrm{for }x\in J_j.
    \end{cases}\]
Then
\begin{align*}
\sum_j \F_L(v_j;\alpha_j,\beta_j) &< \F_{L^c}(u;a,b)
+\sum_j \eta|J_j|\\
&\le  \F_{L^c}(u;a,b) +\eps.
\end{align*}
Hence, since $\int_{\alpha_j}^{\beta_j} |v'|\le \F_L(v_j;\alpha_j,\beta_j)$
and since $v$ agrees with $u$ at the endpoints of the intervals
$J_j$ and outside their union, we see that
$v'$ is integrable and $v$ is an indefinite integral of $v'$.
Consequently, $v\in\AC[a,b]$ and
$\F_L(v;a,b) < \F_{L^c}(u;a,b) +\eps.$

If $x\in[\alpha_j,\beta_j]$, then
  \begin{align*}
    |v(x)-u(x)|& \le  |v(x)-v(\alpha_j)|+|u(x)-u(\alpha_j)|\\
      & \le \F_{L}(v;\alpha_j,\beta_j)+\F_{L^c}(u;\alpha_j,\beta_j)\\
      & <  2\F_{L^c}(u;\alpha_j,\beta_j)+\eta|J_j| \le \eps.
  \end{align*}
So, on recalling that $v(x)=u(x)$ for $x\notin\bigcup_j
[\alpha_j,\beta_j]$, we see that the inequality $|v(x)-u(x)|<\eps$
holds for all $x\in[a,b]$, and the Lemma is proved.

It only remains to prove~(\ref{relax2}). It is enough for us to show
that the statement holds for those $\alpha\in (a,b)$ for which:
\begin{itemize}
\item
there is $\tau>0$ such that
$u$ is affine on $[\alpha-\tau,\alpha+\tau]$;
\item
$\alpha$ is a Lebesgue point of $x\mapsto L(x,u(x),p)$ for every
rational $p$;
\item
$\alpha$ is a Lebesgue point of $x\mapsto L^c(x,u(x),u'(x))$.
\end{itemize}

Let \(\eta>0\). Choose $p_1,p_2\in\R$ and $\lambda_1,\lambda_2\ge 0$
such that $\lambda_1+\lambda_2=1$, $u'(\alpha)=\lambda_1
p_1+\lambda_2 p_2$ and $ \lambda_1 L(\alpha,u(\alpha),p_1)
+\lambda_2 L(\alpha,u(\alpha),p_2) <
L^c(\alpha,u(\alpha),u'(\alpha))+\eta. $ Since $L$ is continuous in
$p$, we can change $p_i$ and $\lambda_i$ slightly so that $p_1$ and
$p_2$ are rational. Let $R=\max(|p_1|,|p_2|)$ and let $C<\infty$ be
an upper bound of $L(x,y,p)$ for $x\in[a,b]$, $|y-u(x)|\le 1$ and
$|p|\le R$.

Choose $\beta\in (\alpha,\alpha+\tau)$ close enough to $\alpha$ so
that
  $$
  \int_\alpha^\beta |L^c(x,u(x),u'(x)) -L^c(\alpha,u(\alpha),u'(\alpha))|\,dx
  <\eta(\beta-\alpha)
  $$
and
  $$
  \int_\alpha^\beta |L(x,u(x),p_i) -L(\alpha,u(\alpha),p_i)|\,dx
  <\eta(\beta-\alpha).
  $$
Since $L(x,y,p)$ is continuous in $y$, there is $0<\delta<1$ so that
the sets
$$S_i=\{x\in(\alpha,\beta): |L(x,y,p_i)-L(x,u(x),p_i)| < \eta
\textrm{ for } |y-u(x)|\le \delta\}\; (i=1,2),$$ each have measure
at least $(1-\eta/C)(\beta-\alpha)$. Choose $m\in\N$  so that
  $$(|p_1|+|p_2|+|u'(\alpha)|)(b-a)/m <\delta.$$
Set $\alpha_j=\alpha+j(\beta-\alpha)/m$,
$\beta_j=\alpha_j+\lambda_1(\beta-\alpha)/m$,
$P_1=\bigcup_{j=0}^{m-1}(\alpha_j,\beta_j)$,
$P_2=\bigcup_{j=0}^{m-1}(\beta_j,\alpha_{j+1})$ and $Q_i=P_i\cap
S_i$, and observe that
  $$|Q_i|\le|P_i|=\lambda_i(\beta-\alpha)\textrm{ and } |P_i\setminus
  Q_i|\le\eta(\beta-\alpha)/C.$$
Define $w\in\AC[\alpha,\beta]$ by $w(\alpha)=u(\alpha)$, $w'=p_1$ on
$P_1$, and $w'(x)=p_2$ on $P_2$. Then $w(\alpha_j)=u(\alpha_j)$ for
$j=1,2,\ldots,m-1$, and $|w(x)-u(x)|<\delta$ for \(x\in[\alpha,\beta
]\). Thus $|L(x,w(x),p_i)-L(x,u(x),p_i)|<\eta$ for $x\in Q_i$, and
so
\begin{align*}
  \lefteqn{\left|\int_{Q_i} L(x,u(x),p_i)\,dx -
  |Q_i|L(\alpha,u(\alpha),p_i)\right|}\\
    &\le \int_{Q_i}|L(x,u(x),p_i)-L(\alpha,u(\alpha),p_i)|\,dx
    \le \eta(\beta-\alpha).
\end{align*}
Hence
\begin{align*}
\int_{P_i} L(x,w(x),p_i)\,dx &\le
\int_{Q_i} L(x,w(x),p_i)\,dx + C|P_i\setminus Q_i|\\
&\le
\int_{Q_i} L(x,u(x),p_i)\,dx + \eta|Q_i|+\eta(\beta-\alpha)\\
&\le
|Q_i|L(\alpha,u(\alpha),p_i)+\eta(\beta-\alpha)+ 2\eta(\beta-\alpha)\\
&\le \lambda_i L(\alpha,u(\alpha),p_i)(\beta-\alpha)+
3\eta(\beta-\alpha).
\end{align*}
Since $\F_{L}(w;\alpha,\beta)$ is the sum of these integrals
over $i=1,2$, we get
\begin{align*}
\F_{L}(w;\alpha,\beta) &\le (\lambda_1 L(\alpha,u(\alpha),p_1)
+\lambda_2 L(\alpha,u(\alpha),p_2))(\beta-\alpha) +6\eta(\beta-\alpha)\\
&\le
L^c(\alpha,u(\alpha),u'(\alpha))(\beta-\alpha) + 7\eta(\beta-\alpha)\\
&\le \F_{L^c}(u;\alpha,\beta) +8\eta (\beta-\alpha),
\end{align*}
and, noting that $w(\beta)=u(\beta)$, we see that the required
statement holds with $u_{\alpha,\beta}=w$, provided that the above
construction was started with $\eta/9$.
\end{proof}

This Lemma clearly implies that
$\F_{L^c}(a,A;b,B)= \F_{L}(a,A;b,B)$ for any $a<b$ and $A,B\in\R$,
and that every relaxed minimizer is
a constrained minimizer. This also implies that every
constrained minimizer is
a relaxed one: if $u\in\AC[a,b]$, $u(a)=A$, $u(b)=B$ is
a constrained minimizer, then
$\F_{L^c}(a,A;b,B)\le \F_{L^c}(u;a,b)\le
\F_{L}(u;a,b)=\F_{L}(a,A;b,B)=\F_{L^c}(a,A;b,B).$
We record this in

\begin{thm}\label{bk6}
If the Lagrangian $L$ is continuous in $(y,p)$ for each \(x\in\R\),
then a function $u\in\AC[a,b]$ is a constrained minimizer for $L$ if
and only if it is a minimizer for $L^c$.
\end{thm}

An immediate consequence of this result and Proposition~\ref{bk2} is
the following theorem.

\begin{thm}\label{generelaxed}
Let $L$ be a Lagrangian that is continuous in $(y,p)$ for each
\(x\in\R\). Let $a<b$ and $A,B$ be given and suppose that there is a
relaxed minimizer $\hat{u}$ such that $\hat{u}(a)=A$,
$\hat{u}(b)=B$,
$$ \liminf_{x\to a}\frac{|\hat{u}(x)-A|}{|x-a|} < \infty \text{ and }
\liminf_{x\to b_-}\frac{|\hat{u}(x)-B|}{|x-b|} < \infty.
$$
Then every generalized minimizer with $u(a)=A$ and $u(b)=B$ is
also a relaxed one.
\end{thm}

The description of the universal singular set of a Lagrangian as the
universal singular set of the convexified Lagrangian also follows
from Theorem~\ref{bk6}.

\begin{cor}\label{univrelaxed}
If \(L\)  is a Lagrangian that is continuous in $(y,p)$ for each
\(x\in\R\), then the universal singular set of $L$
 coincides with the universal singular set of $L^c$.
\end{cor}

\subsection{Tonelli regularity}\label{tr}

In this section, we prove a version of Tonelli's partial regularity
theorem that is valid for functions that are close to minimizers.
The main point of Tonelli's theorem is that, when the slope of a
minimizer $u$ has, between two points $\alpha,\beta$ of its domain,
a certain bound and $\alpha,\beta$ are close enough (depending on
the bound of the slope), then the minimizer is Lipschitz between the
points, and even $|u'|\le C$ on $[\alpha,\beta]$ where $C$ depends
only on the bound and the points. Our idea is that, when the slope
of a given function $u$, between two points $\alpha,\beta$ of its
domain, has a certain bound and $\alpha,\beta$ are close enough
(depending on the bound for the slope), then the measure of the set
$\{x\in[\alpha,\beta]:|u'(x)|> C\}$ (where $C$ depends only on the
bound and the points) should be controlled by the excess.

In the rest of this section, in addition to our usual assumptions on
the Lagrangian $L$ (Borel, bounded from below, locally bounded from
above, and superlinear), we assume that $L$ is locally Lipschitz in
$y$ uniformly for $(x,p)$ in any compact set, that is:
\begin{enumerate}
\item[(L)] For every $R>0$, there is $C\geq 0$ such that
$|L(x,y_1,p)-L(x,y_2,p)|\le C |y_1-y_2|$ whenever\label{test}
$|x|,|y|,|p|\le R$.
\end{enumerate}
 This is weaker than the assumptions in the papers~\cite{tonelli}
 and~\cite{ball_nad}.

\begin{lem}\label{regul.lem}
Let $L$ be a Lagrangian satisfying the Lipschitz
condition~\textup{(L)}. For each $R>0$, there are $M,\delta
>0$ such that if
\begin{enumerate}
  \item $[\alpha,\beta]\subset[a,b]\subset [-R,R]$ with
  $|b-a|<\delta$, and
  \item $u\in\AC[a,b]$ satisfies $|u(x)|\le R$ for $x\in[a,b]$ and
  $$\EE(u;\alpha,\beta) +|u(\beta)-u(\alpha)|+
  \left|\int_{\{x\in[a,b]:|u'(x)|>M\}} u'(x)\,dx\right|\le
  R(\beta-\alpha),
  $$
\end{enumerate}
then
$$
\int_{\{x\in[a,b]:|u'(x)|>M\}} L(x,u(x),u'(x))\,dx \le 2 \EE(u;a,b).
$$
\end{lem}

\begin{proof}
Since adding a positive constant to $L$ does not change the validity
of the Lemma's hypotheses and only strengthens the conclusion, we
may assume without any loss of generality that $L(x,y,p)\ge |p|$ for
all $x,y,p$.

Let \(R>0\) be given.

Fix $C\ge 1$ for which $L(x,y,p)\le C$ whenever $|x|,|y|,|p|\le R$.
Choose $N\ge R+1$ so that
  $$\omega(p)\ge 2(R+C)|p|\text{ for }|p|\ge N,$$
and pick $D\ge C R$ for which $L(x,y,p)\le D$  whenever $|x|\le R$,
$|y|\le N$ and $|p|\le N+2R$.

We can now define the required constants $M$ and $\delta$: first,
choose $M\ge N+2R$ so that $\omega(p)\ge 10\left(1+\frac D R\right)
|p|$ for $|p|\ge M$ and then choose $\delta>0$ so that
  $$|L(x,y_1,p)-L(x,y_2,p)|\le |y_1-y_2|/\delta$$
whenever $|x|\le R$, $|y_1|,|y_2|\le 5R+4RM$ and $|p|\le M$.

Now suppose that intervals $[\alpha,\beta]\subset[a,b]$ and a
function $u\in\AC[a,b]$ satisfy the assumptions of the lemma.
Denoting by $\hat u$ the affine function on $[\alpha,\beta]$ for
which ${\hat u}(\alpha)=u(\alpha)$ and ${\hat u}(\beta)=u(\beta)$,
we infer that
\begin{equation}\label{fr.0}
\F(u;\alpha,\beta)\le \EE(u;\alpha,\beta)+\F(\hat u;\alpha,\beta)\le (R+C)(\beta-\alpha).
\end{equation}
Let $Y=\{x\in[\alpha,\beta]: |u'(x)|>N\}$ and observe that, by the
definition by $N$,
$$
2(R+C)\int_Y |u'(x)|\,dx
\le
\int_\alpha^\beta L(x,u(x),u'(x))\,dx
\le (R+C)(\beta-\alpha).
$$
Hence $|Y|\le (\beta-\alpha)/(2N)\le \tfrac{1}{2}(\beta-\alpha)$ and
\begin{equation}\label{N.ineq}
|\{x\in[\alpha,\beta]: |u'(x)|\le N\}|=|[\alpha,\beta]\setminus
Y|\ge \tfrac{1}{2}(\beta-\alpha).
\end{equation}

Let $Z=\{x\in[a,b]: |u'(x)|>M\}$. For future reference, we note that
\begin{equation}\label{fr.1}
10\left(1+\frac D R\right)\int_Z|u'(x)|\,dx\le \int_Z \omega(u'(x))\,dx
\le \int_Z L(x,u(x),u'(x)) \,dx,
\end{equation}
which implies, in particular, that
\begin{equation}\label{fr.2}
C|Z|\le \frac D R\int_Z|u'(x)|\,dx\le\frac1{10}\int_Z L(x,u(x),u'(x)) \,dx.
\end{equation}
We also notice that the assumption $|u|\le R$ implies that for all
$x\in[a,b]$,
\begin{eqnarray}\label{fr.3}
\left|\int_{\{t\in[a,x]:|u'(t)|>M\}} u'(t)\,dt\right|
&=&
\left|\int_a^x u'(t)\,dt -\int_{\{t\in[a,x]:|u'(t)|\le M\}} u'(t)\,dt\right|\\
&\le& |u(x)-u(a)|+M(x-a)\le 2R(1+M).\nonumber
\end{eqnarray}

We now define a function $v\in\AC[a,b]$ to provide an estimate of
the excess $\EE(u;a,b)$. Combining our hypothesis that $\left|\int_Z
u'(x)\,dx\right|\le R(\beta-\alpha)$ with~\eqref{N.ineq} gives
$$\left|\int_Z u'(x)\,dx\right|
\le 2R|\{x\in[\alpha,\beta]: |u'(x)|\le N\}|.$$
Hence we may find a measurable set
$X\subset\{x\in[\alpha,\beta]: |u'(x)|\le N\}$ with
$$|X|=\frac1{2R}\left|\int_Z u'(x)\,dx\right|.$$
Denoting the sign of $\int_Z u'(x)\,dx$ by $\sigma$, we define
$\phi\colon [a,b]\to\R$ by
$$
  \phi(x)=
    \begin{cases}
      u'(x),    &\textrm{if } x\in Z,\\
      -2\sigma R,&\textrm{if } x\in X,\\
      0,        &\textrm{if } x\in[a,b]\setminus(X\cup Z).
    \end{cases}
$$
Define $v\in\AC[a,b]$ by
  $$v(x)=u(x)-\int_a^x
      \phi(t)\,dt.$$
The  estimates~\eqref{adm.1}--\eqref{adm.4} that follow compare
various energy integrals of $u$ and $v$.

We use two estimates of $|u(x)-v(x)|$ for $x\in[a,b]$: the obvious
one that $|u(x)-v(x)|\le 2\int_Z |u'(t)|\,dt$ and, a consequence
of~\eqref{fr.3}, that $|u(x)-v(x)|\le 4R(M+1)$, which implies
$|v(x)|\le 5R+4RM$, since \(|u(x)|\leq R\).

Our choice of \(\delta \) together with the observations that
$|v'(x)|\le \max(0,N+2R,M)=M$ and $|u(x)-v(x)|\le
2\int_Z|u'(t)|\,dt$ for \(x\in [a,b]\), and~\eqref{fr.1} give
\begin{eqnarray}
\lefteqn{\int_a^b (L(x,v(x),v'(x))-L(x,u(x),v'(x))) \,dx}\qquad\qquad\nonumber\\
&\le& \frac{1}{\delta}\int_a^b |u(x)-v(x)| \,dx
\le \frac{2(b-a)}{\delta}\int_Z|u'(t)|\,dt \nonumber\\
&\le&  2\int_Z|u'(t)|\,dt\le \frac15 \int_Z
L(x,u(x),u'(x))\,dx.\label{adm.1}
\end{eqnarray}

For $x\in X$, we have $|u(x)|\le R$, $|u'(x)|\le N$ and $|v'(x)|\le
N+2R$, so the definitions of the constant $D$ and the set $X$
together with~\eqref{fr.1} give
\begin{eqnarray}
\lefteqn{\int_X (L(x,u(x),v'(x))-L(x,u(x),u'(x))) \,dx}\qquad\qquad\nonumber\\
&\le& 2D|X|\le \frac D R \int_Z|u'(t)|\,dt\le \frac15 \int_Z
L(x,u(x),u'(x))\,dx.\label{adm.2}
\end{eqnarray}

For $x\in Z$, we have $|u(x)|\le R$ and $v'(x)=0$, so
$L(x,u(x),v'(x))=L(x,u(x),0)\le C$. Hence, by \eqref{fr.2},
\begin{eqnarray}
\lefteqn{\int_Z (L(x,u(x),v'(x))-L(x,u(x),u'(x))) \,dx}\qquad\nonumber\\
&\le& C|Z| -\int_Z L(x,u(x),u'(x))\,dx \le -\frac9{10} \int_Z
L(x,u(x),u'(x))\,dx.\label{adm.3}
\end{eqnarray}

Finally, $u'(x)=v'(x)$ for $x\in[a,b]\setminus(X\cup Z)$, and so
\begin{equation}\label{adm.4}
\int_{[a,b]\setminus(X\cup Z)} (L(x,u(x),v'(x))-L(x,u(x),u'(x))) \,dx =0.
\end{equation}

Adding the inequalities \eqref{adm.1}--\eqref{adm.4}, we find
$$
\int_a^b (L(x,v(x),v'(x))-L(x,u(x),u'(x))) \,dx \le
-\frac12 \int_Z L(x,u(x),u'(x))\,dx.
$$
Since $v(a)=u(a)$ and $v(b)=u(b)$, we conclude
\begin{align*}
\int_Z L(x,u(x),u'(x))\,dx
&\le 2\int_a^b (L(x,u(x),u'(x))-L(x,v(x),v'(x))) \,dx\\
&\le 2\EE(u;a,b),
\end{align*}
as required.
\end{proof}

\begin{cor}\label{regul.0}
Let $L$ be a Lagrangian satisfying the Lipschitz
condition~\textup{(L)}. For each $R>0$, there are $M,\delta
>0$ such that if
\begin{enumerate}
  \item $[\alpha,\beta]\subset[a,b]\subset [-R,R]$ with
  $|b-a|<\delta$, and
  \item $u\in\AC[a,b]$ satisfies $|u(x)|\le R$ for $x\in[a,b]$ and
  $$\EE(u;a,b) +|u(\beta)-u(\alpha)|\le R(\beta-\alpha),
  $$
\end{enumerate}
then
$$
\int_{\{x\in[a,b]:|u'(x)|>M\}} L(x,u(x),u'(x))\,dx \le 2 \EE(u;a,b).
$$
\end{cor}

\begin{proof}
Again, we may assume without any loss of generality that
$L(x,y,p)\ge |p|$ for all $x,y,p$.

Let $C>0$ be such that $L(x,y,p)\le C$ whenever $|x|,|y|,|p|\le R$.
We show that the statement holds for $M,\delta$ obtained from
Lemma~\ref{regul.lem} used with $R$ replaced by $R_0:=4(R+C)$. For
this, assume that $[\alpha,\beta]$ and $u\in\AC[a,b]$ are as in the
assumptions of the Corollary.

Call an interval $[c,d]$ \emph{good} if
$[a,b]\supset[c,d]\supset[\alpha,\beta]$ and
$$
\int_{\{x\in[c,d]:|u'(x)|>M\}} L(x,u(x),u'(x))\,dx\le
3R(\beta-\alpha).
$$
Notice that, by Lemma~\ref{regul.lem} (used with $R=R_0$ and
$[a,b]=[c,d]$), if $[c,d]$ is a good interval, then
  $$
    \int_{\{x\in[c,d]:|u'(x)|>M\}}
      L(x,u(x),u'(x))\,dx \le 2 \EE(u;c,d).
  $$
In particular, it suffices to show that $[a,b]$ is good. Since $2
\EE(u;c,d)<3R(\beta-\alpha)$, continuity of the integral implies
that every good interval distinct from $[a,b]$ is contained in a
larger good interval. Hence $[a,b]$ is good, provided we show that
at least one good interval exists. Thus it is enough to show that
$[\alpha,\beta]$ is good.

Let $\hat u$ denote the affine function on $[\alpha,\beta]$ for
which ${\hat u}(\alpha)=u(\alpha)$ and ${\hat u}(\beta)=u(\beta)$.
Then
  $$\int_\alpha^\beta |u'(x)|\,dx \le \F(u;\alpha,\beta)\le
    \EE(u;\alpha,\beta)+\F(\hat{u};\alpha,\beta)\le (R+C)(\beta-\alpha).$$
Hence Lemma~\ref{regul.lem} gives
    $$\int_{\{x\in[\alpha,\beta]:|u'(x)|>M\}} L(x,u(x),u'(x))\,dx \le
        2 \EE(u;\alpha,\beta)\le 3R(\beta-\alpha),$$
as required.
\end{proof}

\begin{cor}\label{singul.1}
Let $L$ be a Lagrangian satisfying the Lipschitz
condition~\textup{(L)}. For every $R,N>0$, there are $M,\delta
>0$ such that if $[a,b]\subset [-R,R]$ with $|b-a|<\delta$, and
if $u\in\AC[a,b]$ is a generalized minimizer for $L$ with $|u|< R$
on $[a,b]$, then either
\begin{itemize}
\item
$|u(\beta)-u(\alpha)|\le M(\beta-\alpha)$ whenever
$a\le\alpha\le\beta\le b$, or
\item
$|u(\beta)-u(\alpha)|\ge N|\beta-\alpha|$ whenever
$a\le\alpha\le\beta\le b$.
\end{itemize}
\end{cor}

\begin{proof}
We again may assume that $L(x,y,p)\ge |p|$ for all $x,y,p$. We show
that the conclusion holds for $M$ and $\delta$ obtained from
Corollary~\ref{regul.0} when $R$ is  replaced by $R_0:=R+N$.

Let $[a,b]$ and $u$ be as in the assumptions and find sequences
$u_k\in\AC[a_k,b_k]$ so that: $u_k\to u$, $\EE(u_k;a_k,b_k)\to 0$,
and $[a_k,b_k]\nearrow [a,b]$.

If the second alternative of the corollary does not hold, then there
are $a\le\alpha < \beta\le b$ for which $|u(\beta)-u(\alpha)|<
N(\beta-\alpha)$. Without loss of generality we can assume that
$a<\alpha<\beta<b$. Then, for sufficiently large $k$,
$[a_k,b_k]\supset[\alpha,\beta]$, $|u_k|\le R$, and
  $$\EE(u_k;a_k,b_k)+|u_k(\beta)-u_k(\alpha)|\le
  R_0(\beta-\alpha).$$
Hence, by Corollary~\ref{regul.0},
\begin{align*}
\int_{\{x\in[a_k,b_k]:|u_k'(x)|>M\}} |u_k'(x)|\,dx&\le
\int_{\{x\in[a_k,b_k]:|u_k'(x)|>M\}} L(x,u_k(x),u_k'(x))\,dx\\
&\le 2 \EE(u_k;a_k,b_k) \to 0,
\end{align*}
showing that the first alternative of the Corollary holds.
\end{proof}

\begin{rem}
Lemma~\ref{regul.lem} and Corollaries~\ref{regul.0}
and~\ref{singul.1}  also hold without change of proof in the
vector-valued case, extending the work of~\cite{CV}. In the
real-valued case more precise information is  apparently available
(see Proposition~\ref{singul.2} following); however it follows
immediately from Corollary~\ref{singul.1}.
\end{rem}

\begin{prp}\label{singul.2}
Let $L$ be a Lagrangian satisfying the Lipschitz
condition~\textup{(L)}. Then for every $R,N>0$, there are $M,\delta
>0$ such that if $[a,b]\subset [-R,R]$ with
$|b-a|<\delta$, and if $u\in\AC[a,b]$ is a generalized minimizer for
$L$ with $|u|< R$ on $[a,b]$, then either
\begin{itemize}
\item
$|u(\beta)-u(\alpha)|\le M(\beta-\alpha)$ whenever
$a\le\alpha<\beta\le b$, or
\item
$u(\beta)-u(\alpha)\ge N(\beta-\alpha)$ whenever
$a\le\alpha<\beta\le b$, or
\item
$u(\beta)-u(\alpha)\le -N(\beta-\alpha)$ whenever
$a\le\alpha<\beta\le b$.
\end{itemize}
\end{prp}

\begin{proof}
Let $M$ and $\delta$ be as given by Corollary~\ref{singul.1} for
\(R\). Corollary~\ref{singul.1} implies that if the  first
alternative does not hold, then $|u(\beta)-u(\alpha)|\ge
N(\beta-\alpha)$ whenever $a\le\alpha<\beta\le b$. In this case it
is enough to show that, for $a\le\alpha<\beta\le b$,
$u(\beta)-u(\alpha)$ is either always positive or always negative.
But this is obvious, since otherwise there are
$a\le\alpha<\beta\le b$ such that $u(\beta)-u(\alpha)=0$ in which
case Corollary~\ref{singul.1} implies that the first alternative
occurs.
\end{proof}

Another way of presenting our results is closer to the usual
formulation of Tonelli's partial regularity theorem. It follows
directly from the fact that $u$ has a finite derivative almost
everywhere and the previous Proposition.

\begin{prp}\label{singul.7}
Let $L$ be a Lagrangian satisfying the Lipschitz condition {\rm
(L)}. Then for every generalized minimizer $u$ on $[a,b]$, there
are disjoint closed Lebesgue null sets $E_+,E_-\subset[a,b]$ such
that
\begin{itemize}
\item $u$ is locally Lipschitz on $[a,b]\setminus (E_+\cup E_-)$;
\item $\lim_{s\not= t,\, \max (d(s,E_+),d(t,E_+),\, |t-s|)\to 0} (u(t)-u(s))/(t-s) = \infty$;
\item $\lim_{s\not= t,\, \max (d(s,E_-),d(t,E_-),\, |t-s|)\to 0} (u(t)-u(s))/(t-s) = -\infty$.
\end{itemize}
\end{prp}

Many standard variants of Tonelli's regularity results may be
obtained, under appropriate smoothness and strict convexity
assumptions, by deducing the Euler-Lagrange equation on the
intervals where the minimizer is Lipschitz. Since this is a
straightforward use of known methods, we do not do this here.

\section{The structure of universal singular sets}\label{su}

This section is devoted to the study of the size of the intersection
of universal singular sets with rectifiable curves. Our aim is to
show that universal singular sets intersect many curves in a set of
zero length.

Recall that the universal singular set of \(L\) is defined to be all
points \( (x_0,y_0)\in\R^2\) for which there are $a\leq  x_0\leq b$
with \(a<b\) and a  \(u\in \AC[a,b]\) that is a generalized
minimizer for \(L\) with $u(x_0)=y_0$ and $|u'(x_0)|=\infty$.
However, for continuous Lagrangians that are convex in $p$ our
argument does not use generalized minimizers.

It turns out that, although the universal singular set of any
Lagrangian meets the graph of every absolutely continuous function
in a set of linear measure zero, the situation is considerably
more delicate when it comes to curves that may have vertical
tangents: vertical lines also meet the universal singular set in a
null set, but, as shown in section~\ref{ll}, some rectifiable
curves may actually meet it in a set of positive linear measure.

In this section, we consider the following problems:
\begin{enumerate}
    \item Which curves have the property that they meet the
    universal
    singular set of each Lagrangian with a given superlinear growth in
    a null set?\label{q1}
    \item Which curves have the property that they
    meet the universal singular set of every Lagrangian in a null
    set?\label{q2}
    \item Which curves have the property that they meet a set
    contained in the universal singular set of Lagrangians with
    arbitrary superlinear growth in a null set?\label{q3}
    \item For which
    Lagrangians is the universal singular set of the first
    category?\label{q4}
\end{enumerate}

We now describe our results. The examples of the following section
show that these results are close to a complete picture.

Our answer to problem~(\ref{q1}) is given by the following theorem.
\begin{thm}\label{thm.omega-pu}
    Let $\omega\colon\R\to\R$ be even, convex and superlinear. Suppose
    that an absolutely continuous curve $\gamma(t)=(x(t),y(t))\colon
    [a,b]\to\R^2$ is such that for almost all $t\in [a,b]$, either
    \begin{equation}\label{thm.omega-pu.1}
    \limsup_{s\to t} \left|\frac{y(s)-y(t)}{x(s)-x(t)}\right|<\infty
    \end{equation}
    or
    \begin{equation}\label{thm.omega-pu.2}
    \liminf_{s\to t} |x(s)-x(t)|\;
    \omega\left(\frac{y(s)-y(t)}{x(s)-x(t)}\right)>0.
    \end{equation}
    Then \( \{\gamma(t):t\in[a,b]\}\) meets the universal singular set
    of any Lagrangian~\(L\) for which $L(x,y,p)\ge\omega(p)$ in a set
    of linear measure zero.
    (When $x(s)-x(t)=0$, we take
    $$\left|\frac{y(s)-y(t)}{x(s)-x(t)}\right|\mbox{ to be zero and }
      |x(s)-x(t)|\omega\left(\frac{y(s)-y(t)}{x(s)-x(t)}\right)\mbox{ to
      be }\infty.\mbox{)}$$
\end{thm}

Before proving this Theorem, we show how it is used to answer
questions~(\ref{q2})--(\ref{q4}). Since absolutely continuous
functions (considered as curves $x\to (x,f(x))$)
satisfy~(\ref{thm.omega-pu.1}), and vertical lines
satisfy~(\ref{thm.omega-pu.2}), we have an answer to
question~(\ref{q2}).

\begin{thm}\label{thm.near-pu}
Graphs of absolutely continuous functions and vertical lines meet the
universal singular set of any Lagrangian in a set of linear measure zero.
\end{thm}

The answers to  problems~(\ref{q3}) and~(\ref{q4}) also follow from
Theorem~\ref{thm.omega-pu}, but a little more work is needed.

\begin{thm}\label{thm.full-pu}
Suppose that $E\subset\R^2$ is such that for any superlinear
$\omega$ there is a Lagrangian $L(x,y,p)\ge\omega(p)$ whose universal
singular set contains $E$. Then $E$ is purely unrectifiable.
\end{thm}

\begin{proof}
It suffices to show that $|\{t\in(a,b): \gamma(t)\in E\}|=0$ for any
injective $C^1$ curve $\gamma(t)=(x(t),y(t))\colon [a,b]\to\R^2$ for
which $\gamma'(t)\ne 0$ for all $t\in[a,b]$. This follows from
Theorem~\ref{thm.omega-pu} provided we find a superlinear function
$\omega$ for which
\begin{equation}\label{fp.p1}
\liminf_{s\to t}|x(s)-x(t)|\
\omega\left(\frac{|y(s)-y(t)|}{|x(s)-x(t)|}\right)>0\text{
whenever }x'(t)=0.
\end{equation}

Let $\delta_n\searrow 0$ be a strictly decreasing sequence such
that for \(s,t\in [a,b]\) with \(s\not=t\),
    $$|y(s)-y(t)|>(n+1)|x(s)-x(t)|\text{ whenever }x'(t)=0\text{
    and }|x(s)-x(t)|<\delta_n.$$
 Such a sequence exists, since otherwise we
could find, for some $n$, sequences $s_k,t_k\in[a,b]$ with
$$x'(t_k)=0,\, |y(s_k)-y(t_k)|\le (n+1)|x(s_k)-x(t_k)|\text{ and }
    |x(s_k)-x(t_k)|\to 0.$$
 But then, by passing to a subsequence, we could assume that $s_k\to
s$ and $t_k\to t$ with
    $$x'(t)=0,\, |y(s)-y(t)|\le (n+1)|x(s)-x(t)|\text{ and }
        x(s)=x(t).$$
 But this means $y(s)=y(t)$, and so, since $\gamma$ is injective, $s=t$. But then
$|y'(t)|>0=x'(t)$, and so $|y(s_k)-y(t_k)|>(n+1)|x(s_k)-x(t_k)|$ for
$k$ large enough --- a contradiction.

Let $\kappa\colon(0,\infty)\to(0,\infty)$ be a continuous
decreasing function for which $\kappa(\delta_n)=n$. Let
$\kappa^{-1}$ be the inverse of $\kappa$ and let $\omega$ be any
convex superlinear function satisfying $\omega(p)=\omega(|p|)\ge
1/\kappa^{-1}(p)$. If $s,t\in(a,b)$ and $\delta_{n+1}\le
|x(s)-x(t)|<\delta_n$, then
    $$\frac{|y(s)-y(t)|}{|x(s)-x(t)|}\ge n+1=\kappa(\delta_{n+1})\ge\kappa(|x(s)-x(t)|).$$
Hence
    $$|x(s)-x(t)| \omega\left(\frac{|y(s)-y(t)|}{|x(s)-x(t)|}\right)
    \ge|x(s)-x(t)|\omega(\kappa(|x(s)-x(t)|))\ge 1, $$
 for every
$s\in[a,b]$ for which $|x(s)-x(t)|<\delta_1$, and so, in
particular, for $|s-t|$ small enough. It follows
that~\eqref{fp.p1} holds and we are done.
\end{proof}

Our answer to problem~(\ref{q4}) is given in the following Theorem
where we show that,  assuming the Lagrangian satisfies the Lipschitz
condition~\textup{(L)}, the universal singular set is of the first
Baire category. Our proof
 is based on the regularity results of
Section~\ref{tr} and the  `almost' pure unrectifiability of
universal singular sets described in Theorem~\ref{thm.omega-pu}. We
show in Section~\ref{ll} that some additional assumption on the
Lagrangian is necessary.

\begin{thm}\label{f-sigma}
If $L$ is a Lagrangian that satisfies the Lipschitz
condition~\textup{(L)}, then the universal singular set of $L$ is a
countable union of closed sets. In particular, it is a first
category set.
\end{thm}

\begin{proof}
Let $L$ be a Lagrangian that satisfies the Lipschitz
condition~\textup{(L)}. Let $S_k$ be the set of all $P=(a,A)\in\R^2$
for which there is a generalized minimizer $u^P$ on $[a,a+\tfrac1k]$
with $|u^P|\le k$ and $u'(a)=+\infty$. By
Proposition~\ref{singul.2}, for every $N>0$, there is
$0<\delta_N<\tfrac1k$ such that
    $$u^P(y)-u^P(x)\ge N(y-x)\text{ whenever }
    a\le x<y\le a+\delta_N.$$
 So if $P_l \in S_k$ converge to
some $P\in\R^2$, then the corresponding functions $u^{P_l}$ converge
(up to a subsequence) to some generalized minimizer $u$ on
$[a,a+1/k]$ for which $|u|\le k$, and
  $$u(y)-u(x)\ge N(y-x)\mbox{ whenever }a\le x<y\le
    a+\delta_N.$$
Hence $u'(a)=+\infty$, and so $P\in S_k$ implying that $S_k$ is a
closed set.

The universal singular set of \(L\) is a union of
$\bigcup_{k=1}^\infty S_k$ together with three other sets obtained
by symmetrical constructions. Since each of these sets is a
countable unions of closed sets, the universal singular set is a
countable union of closed sets.

Moreover, Theorem~\ref{thm.omega-pu} implies that the interior of
the universal singular set of \(L\) is empty. Hence the universal
singular set of \(L\) is a first category set, since a countable
union of closed sets in the plane that is not of the first category
has nonempty interior.
\end{proof}

The key to the proof of Theorem~\ref{thm.omega-pu}, and hence these
results, is in understanding that, for fixed \( (a,A)\in\R^2\), the
functional $(x,X)\mapsto\F(a,A;x,X)$ increases (or decreases)
steeply in many directions from $(b,B)$ whenever \( (b,B)\) is in
the universal singular set. This result is of independent interest.

\begin{lem}\label{phi3a}
Let $\omega\colon\R\to\R$ be convex and superlinear, and  suppose
that $L$ is a Lagrangian satisfying $L(x,y,p)\ge\omega(p)$ for all
\(x,y,p\).

If $a,A,b,B\in\R$ with $a<b$, and  $u\in \AC[a,b]$ satisfies
$u(a)=A$, $u(b)=B$ and $u'_-(b)=\infty$, then for all \(C,D>0\)
  $$
    \lim_{T\ni(x,X)\to (b,B)}
    \frac{\E(a,A;x,X)-\F(u;a,b)}{\|(x,X)-(b,B)\|} =-\infty,
  $$
  where
  \begin{align*}
  T&:=\left\{(x,X): x<b,\,X<B\text{ and } (b-x)\omega\left(\frac{B-X}{b-x}\right) > D\right\}\\
  &\qquad {}\cup \{(x,X): x\ge b\text{ and } X-B< C(x-b) \}.
  \end{align*}
\end{lem}

\begin{proof}
Since $L$ is bounded from below, $\F(u;a,b)>-\infty$ and so, the
case $\F(u;a,b)=\infty$ being trivial, we may assume that
$\F(u;a,b)$ is finite.

Fix  $0<K<\infty$; it is enough to show that for any \( (x,X)\in T\)
that is sufficiently close to \( (b,B)\), \(\E(a,A;x,X)<\F(u; a,b)-
K\|(x-b,X-B)\|\).

Let $M=\sup_{x\in[a,b+1], |y|\le |B|+1} |L(x,y,2C)|$. Choose $p_0 >
2C+1$ so that $\omega $ is increasing and $\omega(p)
> (K+2M(C+1))(1+1/C)p$  on $[p_0,\infty)$.
Choose $\alpha\in (\max(a,b-1) ,b)$ sufficiently close to $b$ so
that
\begin{enumerate}
\item\label{ph.1.3} $(b-\alpha)\omega(p)<D$ for $0\le p\le p_0$;
\item\label{ph.1.1} $\F(u;s,t)<D$ whenever $\alpha\le s<t\le b$;
\item\label{ph.1.4} $|u(x)-B|<1$ for $\alpha\le x\le b$;
\item\label{ph.1.2} $u(x) < B + p_0(x-b)$ for $\alpha\le x<b$.
\end{enumerate}

Since~\iref{ph.1.2} with $x=\alpha$ implies $B
+2C(\alpha-b)>u(\alpha)$, there is $0<\delta<b-\alpha$ so that
  $$X_0+2C(\alpha-x_0) > u(\alpha)\text{ whenever }\|(x_0,X_0)-(b,B)\|<\delta.$$

Fix $(x_0,X_0)\in T$ for which $\|(x_0-b,X_0-B)\|<\delta$. Note that
  $$B-X_0 > C(b-x_0);$$
for, if $x_0<b$, then $(b-\alpha)\omega((B-X_0)/(b-x_0)) > D$,
so~\iref{ph.1.3} implies that $B-X_0 > p_0(b-x_0)$, and if $x_0\ge
b$, then the definition of $T$ gives directly that $B-X_0 >
C(b-x_0)$.

Let $\beta=\min(x_0,b)$. We show that
\begin{equation}
X_0+2C(\beta-x_0) < u(\beta).\label{eqn-left}
\end{equation}

If $x_0\ge b$, then \(\beta=b\), \(u(\beta)=B\),
and~(\ref{eqn-left}) follows from $X_0-B < C(x_0-b)<2C(x_0-b)$.

If $x_0<b$, inequality~(\ref{eqn-left}) is just $X_0<u(x_0)$.
Assuming, for a contradiction, that $X_0\ge u(x_0)$, we use
$(B-u(x_0))/(b-x_0)> p_0$, the monotonicity of $\omega$ on
$[p_0,\infty)$ and that \( (x_0,X_0)\in T\) to infer
  $$(b-x_0)\omega((B-u(x_0))/(b-x_0))\ge
  (b-x_0)\omega((B-X_0)/(b-x_0))
  > D.$$
Jensen's inequality then gives
  \begin{align*}
    D> \F(u;x_0,b) &\ge \int_{x_0}^{b}\omega(u'(x))\,dx \\
      &\ge (b-x_0)\omega\left(\frac{u(b)-u(x_0)}{b-x_0}\right)\\
      & \to (b-x_0)\omega\left(\frac{B-u(x_0)}{b-x_0}\right) > D\mbox{ --- a
    contradiction.}
  \end{align*}

Since $X_0+ 2C(\alpha-x_0) > u(\alpha)$ and $X_0+ 2C(\beta-x_0) <
u(\beta)$, there is $\tau\in(\alpha,\beta)\subset (\alpha, b)$ such
that $X_0+ 2C(\tau-x_0) = u(\tau)$.

We use two estimates of $\|(x_0-b,X_0-B)\|$. First, by rearranging
  $$X_0- 2C(b-\tau)-2C(x_0-b) = u(\tau) < B-(2C+1)(b-\tau),$$
we find $b-\tau < B-X_0+ 2C(b-x_0) \le (2C+1)\|(x_0-b,X_0-B)\|.$
Hence
\begin{equation}\label{eph.1.1}
x_0-\tau\le (x_0-b)+ (b-\tau) < 2(C+1)\|(x_0-b,X_0-B)\|.
\end{equation}
For the second estimate, we recall that $\|(x_0-b,X_0-B)\|<\delta$
and $(x_0,X_0)\in T$, and so $X_0-B < C(x_0-b)$. We also notice that
  $$B-u(\tau)=B-X_0+2C(x_0-\tau)\ge B-X_0.$$
Thus if $x_0\le b$, then  $0\le C(b-x_0) < B-X_0$, and so
  $$\|(x_0-b,X_0-B)\|\le (b-x_0)+(B-X_0) < (1+1/C)(B-u(\tau)).$$
Whereas if $x_0>b$, we use
  $$B-u(\tau)=B-X_0+2C(x_0-\tau)\ge B-X_0+2C(x_0-b)> C(x_0-b)$$
to infer that $|X_0-B|\le \max(B-u(\tau),C(x_0-b))=B-u(\tau)$ and
so, in either case,
\begin{equation}\label{eph.1.2}
\|(x_0-b,X_0-B)\|\le (x_0-b)+|X_0-B| < (1+1/C)(u(b)-u(\tau)).
\end{equation}

Let
$$v(x)=\begin{cases}
X_0+ 2C(x-x_0) & \textrm{on } [\tau,x_0],\\
u(x) & \textrm{on } [a,\tau].
\end{cases}$$
Noting that $[\tau,x_0]\subset(b-1,b+1)$ and $|v(x)|\le |B|+1$ on
$[\tau,x_0]$, we infer from~\eqref{eph.1.1} that
  $$\F(v; \tau,x_0)\le M(x_0-\tau) < 2M(C+1)\|(x_0-b,X_0-B)\|.$$
However Jensen's inequality and~\eqref{eph.1.2} imply
\begin{align*}
\F(u;\tau,b)\ge  \int_{\tau}^{b} \omega(u'(x))\,dx &\ge
(b-\tau)\omega\left(\frac{u(b)-u(\tau)}{b-\tau}\right)\\
&> (K+2M(C+1))\|(x_0-b,X_0-B)\|.
\end{align*}
Hence
\begin{align*}
  \F(v; a,x_0) &\le  \F(u; a,b)-\F(u; \tau,b)+\F(v; \tau,x_0)\\
    &< \F(u; a,b)- K\|(x_0-b,X_0-B)\|,
\end{align*}
giving $\E(a,A;x_0,X_0)<\F(u; a,b)- K\|(x_0-b,X_0-B)\|$,
as required.
\end{proof}

In order to prove Theorem~\ref{thm.omega-pu}, we use the following
well-known result showing that certain sets of reals always have
linear measure zero,. (See~\cite{saks} for many ramifications.)

\begin{lem}\label{inf-der}
The set of points at which a function $f\colon [a,b]\to\R$ has an
infinite one-sided derivative has linear measure zero.
\end{lem}

\begin{proof}
It suffices to consider the set $E=\{x\in (a,b): f'_+(x)=+\infty\}.$
Let $E_n=\{x\in E: f(x)\le f(y)\text{ for }
x<y<x+\frac{1}{n}(b-a)\}$ and
$E_{n,k}=E_n\cap[a+(k-1)(b-a)/n,a+k(b-a)/n).$ Then
$E=\bigcup_{n=1}^\infty\bigcup_{k=1}^n E_{n,k}$, and so it is enough
to show that each $E_{n,k}$ has measure zero.

Fix $n,k\in \N$ and observe that $f(x)\le f(y)$ whenever $x,y\in
E_{n,k}$ and $x<y$. Hence there is a nondecreasing function $g$ that
agrees with $f$ on $E_{n,k}$. If $x\in E_{n,k}$ is not a
right-isolated point of $E_{n,k}$ (that is, $(x,x+\delta)\cap
E_{n,k}\ne\emptyset$ for every $\delta>0$), then  $g$ does not have
a finite derivative at \(x\). But $g$, being monotonic,
 has a finite derivative almost everywhere, and so $E_{n,k}$  has
measure zero since it can have at most countably many right-isolated
points.
\end{proof}

\begin{proof}[Proof of Theorem \ref{thm.omega-pu}]
Given $(a,A)\in\R^2$, denote by $S_{a,A}$ the set of $(b,B)\in\R^2$
for which $a<b$ and there is $u\in\AC[a,b]$ that is a constrained
minimizer  for $L$ with $u(a)=a$, $u(b)=B$ and $u'_-(b)=+\infty$. It
is enough to show that the union of the sets, $S:=\bigcup_{a,A}
S_{a,A}$, is null on the curves described in the Theorem; the whole
universal singular set is a union of this set and three other sets
that are given by similar definitions, and to which symmetrical
arguments apply.

In fact we can assume that the union defining $S$ is only taken over
rational $a$ and $A$. For given $(b,B)\in S_{a,A}$ and a
corresponding constrained minimizer $u\in\AC[a,b]$, we pick a
rational $\hat a\in (a,b)$ and a rational $\hat A < u(\hat a)$, use
Corollary~\ref{constrained-prop} to infer that $u$ is a constrained
minimizer on $[\hat a,b]$, use Theorem~\ref{exi} to find a
constrained minimizer $v$ on $[\hat a,b]$ with $v(\hat a)=\hat A)$
and $v(b)=B$, and use Proposition~\ref{gm-min} to infer that
$\min(u,v)$ is a constrained minimizer, and conclude that $(b,B)\in
S_{\hat a,\hat A}$.

Hence it is enough to show that $S_{a,A}$ is a null set for all
$a,A\in\R$.

Let $\gamma(t)=(x(t),y(t))$ be an absolutely continuous curve
defined on an interval $[a,b]$ and having the properties
\eqref{thm.omega-pu.1} and \eqref{thm.omega-pu.2} from the Theorem.
Since $\gamma$ maps both the set of points at which it is
non-differentiable and the set of points at which it has derivative
zero to a set of linear measure zero, we restrict our attention to
the set $E$ of those $t\in(\alpha,\beta)$ at which $\gamma$ is
differentiable, has nonzero derivative and has one of the properties
from the Theorem. Then subsets of $E$ have Lebesgue measure zero if
and only if their image under $\gamma$ has linear measure zero; so
it is enough to show that $F:=\{t\in E: \gamma(t)\in S_{a,A}\}$ has
Lebesgue measure zero. We do this by showing that the right
derivative of the function $f(t):=\F(a,A;x(t),y(t))$ is infinite at
every point of $F$.

So let $t\in F$. If~\eqref{thm.omega-pu.1} holds, choose \(C\) so
that
$$\limsup_{s\to t} \left|(y(s)-y(t))/(x(s)-x(t))\right|<C<\infty$$
and $D>0$ arbitrarily. If~\eqref{thm.omega-pu.2} holds, choose \(D\)
so that
$$\liminf_{s\to t} |x(s)-x(t)|
\omega((y(s)-y(t))/(x(s)-x(t)))>D>0$$
and $C<\infty$ arbitrarily.
Let $(b,B)=(x(t),y(t))$ and
\begin{align*}
T&:=\{(x,X): x<b,\,X<B,\, (b-x)\omega((B-X)/(b-x)) > D\}\\
& \qquad {}\cup \{(x,X): x\ge b,\, X-B< C(x-b) \}.
\end{align*}
Then Lemma~\ref{phi3a} applied with any constrained minimizer $u$
witnessing that $(b,B)\in S_{a,A}$, and the fact that
$\F(u;a,b)=\E(a,A;b,B)$, gives
\begin{equation}\label{pf.omega-pu.1}
\lim_{T\ni(x,X)\to (b,B)}
\frac{\E(a,A;x,X)-\E(a,A;b,B)}{\|(x,X)-(b,B)\|} =-\infty.
\end{equation}

Since $\gamma'(t)\ne 0$, we have $0<\lim_{s\to t}
\|((x(s),y(s))-(b,B)\|/|s-t|<\infty$ and one of the following cases
must occur:
\begin{itemize}
\item
If $x'(t)> 0$, then~\eqref{thm.omega-pu.1} holds. Hence for $s>t$
sufficiently close to $t$ we have $x(s)-B< C(y(s)-b)$, and so
$(x(s),y(s))\in T$. Thus~\eqref{pf.omega-pu.1} implies that
$f'_+(t)=-\infty$.
\item
If $x'(t)< 0$, then~\eqref{thm.omega-pu.1} holds. Hence
$(x(s),y(s))\in T$ for $s<t$ sufficiently close to $t$.
Thus~\eqref{pf.omega-pu.1} implies that $f'_-(t)=\infty$.
\item
If $x'(t)= 0$ and $y'(t)>0$, then~\eqref{thm.omega-pu.2} holds, and
we see that $(x(s),y(s))\in T$ for $s>t$ sufficiently close to $t$.
Hence~\eqref{pf.omega-pu.1} implies that $f'_+(t)=-\infty$.
\item
If $x'(t)= 0$ and $y'(t)<0$, then~\eqref{thm.omega-pu.2} holds, and
we see that $(x(s),y(s))\in T$ for $s>t$ sufficiently close to $t$.
Hence~\eqref{pf.omega-pu.1} implies that $f'_-(t)=\infty$.
\end{itemize}

We see that in each case,  $f$ has an infinite one-sided derivative
at $t$. Hence Lemma~\ref{inf-der} shows that $F$ has measure zero,
and the Theorem is proved.
\end{proof}

\section{Lagrangians with large universal singular sets}\label{ll}

In this section we show that the results of
Theorems~\ref{thm.omega-pu}--\ref{f-sigma} are close to being
optimal.

For Theorems~\ref{thm.omega-pu}--\ref{thm.full-pu}, we give examples
of smooth Lagrangians satisfying classical conditions (including
convexity in $p$) for which the universal singular set is as large
as possible. In fact, the Lagrangians in Theorems~\ref{thm.pu}
and~\ref{thm.non-pu} have the following special form: we assume that
we are given a strictly convex superlinear function $\omega\in
C^\infty(\R)$ for which \(\omega (0)=0\), and we construct
Lagrangians $L$ for which
  \begin{enumerate}
    \item[$(\star)$] $L(x,y,p)=\omega(p)+F(x,y,p)$ where $F$ satisfies:
    \begin{enumerate}
      \item[$(\star_1)$] $F\in C^\infty(\R^3)$;
      \item[$(\star_2)$] $F\ge 0$ and for all \(x,y\in\R\), $F(x,y,0)=0$;
      \item[$(\star_3)$] $p\mapsto F(x,y,p)$ is convex for each fixed \( (x,y)\).
    \end{enumerate}
  \end{enumerate}
Notice that for such Lagrangians, the classical existence theorems
hold, and the universal singular set corresponds with that defined
by Ball and Nadirashvili.

The main result of this section is given in the following theorem.
\begin{thm}\label{thm.pu}
  Fix a strictly convex superlinear function $\omega\in C^\infty(\R)$
  for which $\omega(p)\ge\omega(0)=0$, and let $S\subset\R^2$ be a purely
  unrectifiable compact set. Then there is a Lagrangian
  satisfying~$(\star)$ whose universal singular set contains $S$.
\end{thm}

Recall that purely unrectifiable compact subsets of $\R^2$ may have
Hausdorff dimension two and may contain non-trivial continua; so, in
spite of Theorem~\ref{thm.near-pu}, universal singular sets may be
rather large.

We complement this result by a more particular example showing that,
even when one restricts to compact sets, Theorem~\ref{thm.pu} does
not provide a complete answer.

\begin{thm}\label{thm.non-pu}
  Fix a strictly convex superlinear function $\omega\in C^\infty(\R)$
  for which $\omega(p)\ge\omega(0)=0$. Then there is a rectifiable compact set
  $S\subset\R^2$ of positive linear measure that is contained in the
  universal singular set of some Lagrangian satisfying~$(\star)$.
\end{thm}

Unlike the measure zero result of Sych\"ev, our generalisation of
the first category result of Ball and Nadirashvili
(Theorem~\ref{f-sigma}) is shown only under (mild) additional
smoothness assumptions on the Lagrangian. The following result shows
that this is necessary.

\begin{thm}\label{resid}
Fix a
 superlinear function  \(\omega\colon \R\to [0,\infty)\)
for which \( \omega (0)=0\). Then there is a continuous Lagrangian
$L$ with \(L (x,y,p)\geq \omega (p)\) for \((x,y,p)\in\R^3\) and
whose universal singular set is residual in $\R^2$.
\end{thm}

We start by describing the general ideas behind our constructions.
To construct Lagrangians with a large singular set $S$, and with a
given superlinearity $\omega$, we employ the idea of calibrations.
Basically, we prescribe a field of minimizers (better: functions
that should be minimizers in the future) that have infinite
derivative when passing through the points of $S$. These minimizers
are given by the equation $u'=\psi(x,u)$ for a suitable function
$\psi$ that we expect to have singular behaviour whenever $(x,u)\in
S$. We also choose (at this stage completely independently) the
potential of the energy of our field of minimizers; that is, a
function $\Phi\colon\R^2\to\R$. Since we want  each of our (future)
minimizers \(u\) to satisfy $\int_a^b L(x,u,u')\,dx=
\Phi(b,u(b))-\Phi(a,u(a))$, and so \(L\) is a null Lagrangian for
the minimizers \(u\), we have to define  $L(x,y,p) = \Phi_x(x,y)+
p\Phi_y(x,y)$ for $p=\psi(x,y)$. The superlinearity condition means
that we  require
\begin{equation}\label{intr.1.0}
\Phi_x(x,y)+ \psi(x,y)\Phi_y(x,y)\ge\omega(\psi).
 \end{equation}
The calibration argument, which is formally given by
\begin{align*}
  \int_a^b L(x,u,u')\,dx &\ge \int_a^b \Phi_x(x,u(x))+
  \Phi_y(x,u(x))u'(x)\,dx\\
    &=\Phi(b,u(b))-\Phi(a,u(a)),
\end{align*}
leads to the second requirement, namely
\begin{equation}\label{intr.1.1}
L(x,y,p)\ge \Phi_x(x,y)+ p\Phi_y(x,y)\mbox{ whenever }(x,y)\notin S.
\end{equation}
We manage to avoid the set $S$ here because the pure
unrectifiability of $S$ should imply that, along trajectories, the
inequality \eqref{intr.1.1} holds almost everywhere. However the
fact that $\Phi$ cannot behave regularly at the points of $S$
returns to haunt us. The easier difficulty, that~\eqref{intr.1.1}
may have no continuous solution $L$, is avoided rather simply by
requiring that close to $S$, $\Phi_x$ is negative and $|\Phi_y|$ is
much smaller that $|\Phi_x|$; hence the function on the right hand
side of~\eqref{intr.1.1} is locally bounded from above. The harder
problem is that, for some absolutely continuous $u(x)$, the function
$\Phi(x,u(x))$ may fail to be absolutely continuous. To handle this,
we use the pure unrectifiability requirement on $S$ to construct
$\Phi$ in such a way that, at least for increasing $u(x)$ (to which
the problem may be reduced by a simple trick), the composition
$\Phi(x,u(x))$ maps null sets to null sets. Other properties of
\(\Phi(x,u(x))\) and classical real analysis then imply that it is
in fact absolutely continuous, justifying the use of the formal
calibration argument indicated above.

Before embarking on the technical details which, unfortunately,
involve a little more than the above basic description, we should
comment on~\eqref{intr.1.0}.  Since our main discussion did not
involve any condition on $\psi$, we are free to choose it subject
only to~\eqref{intr.1.0}. It is therefore more natural to begin by
defining $\Phi$ satisfying all the requirements alluded to above and
then choose a suitable $\psi$. Noting that the first part of our
argument leads to $\Phi_x<0$, we impose $\Phi_y>0$ (the signs are,
of course, arbitrary; these come from the requirement that the
relevant minimizers be increasing) and observe that the increase of
energy proved in Section~\ref{su} shows that $\psi>-\Phi_x/\Phi_y$
(at least close to the points of the universal singular set, and it
should be substantially bigger there). Hence we decide to take
$\psi=-2\Phi_x/\Phi_y$. This transforms~\eqref{intr.1.0} into
$-\Phi_x\ge\omega(-2\Phi_x/\Phi_y)$, which is easy to achieve. In
reality, technical points force us to require similar but much
stronger inequalities. We therefore begin our argument by giving the
details of the requirements on $\Phi$ leading to the construction of
the required Lagrangian in Lemma~\ref{out}, and only then describe
the particular constructions of $\Phi$ giving the proofs of
Theorems~\ref{thm.pu} and~\ref{thm.non-pu}, respectively.

Before giving the proofs of Theorems~\ref{thm.pu}
and~\ref{thm.non-pu}, we record a couple of results that are used in
the constructions of the required Lagrangians.

The following simple lemma tells us how to smooth the corners of a
particular piecewise-affine function.
\begin{lem}\label{corner}
There is a $C^\infty$ function
$\gamma\colon\{(p,a,b)\in\R^3: b >0\}
\to\R$
such that:
\begin{in.enumerate}\alph
\item\label{corner.1} $p\mapsto \gamma(p,a,b)$ is convex;
\item\label{corner.2} $\gamma(p,a,b)=0$ for $p\le a-1$;
\item\label{corner.3} $\gamma(p,a,b)=b(p-a)$ for $p\ge a+1$;
\item\label{corner.4} $\gamma(p,a,b)\ge \max(0, b(p-a))$.
\end{in.enumerate}
\end{lem}

\begin{proof}
Let $\eta\colon \R\to[0,1]$ be a non-decreasing $C^\infty$ function
such that $\eta(p)=0$ for $p\le -1$,
$\eta(p)=1$ for $p\ge 1$, and $\int_{-1}^1\eta(t)dt=1$.
Define
$\gamma(p,a,b)=b\int_{-\infty}^{p-a} \eta(t)dt.$
Then \iref{corner.1} holds since
$\partial \gamma/\partial p$ is non-decreasing,
\iref{corner.2} is obvious
since for $p\le a-1$ the integrand vanishes,
for \iref{corner.3} we have
$b\int_{-\infty}^{p-a} \eta(t)dt=b\int_{-1}^1\eta(t)dt
+b\int_{1}^{p-a}\eta(t)dt=b(p-a)$,
and \iref{corner.4} follows from the previous statements.
\end{proof}

The following lemma gives sufficient conditions to assert
the existence of a Lagrangian satisfying~\( (\star)\) with a
given compact set inside its universal singular set, and plays
a central role in the proofs of Theorems~\ref{thm.pu} and~\ref{thm.non-pu}.

To simplify notation in what follows, $u$ always denotes a
real-valued function defined on an interval and $U$ denotes the
corresponding function from the same interval to the plane given by
$U(x)=(x,u(x))$.

\begin{lem}\label{out}
Fix a strictly convex superlinear function $\omega\in C^\infty(\R)$
for which $\omega(p)\ge\omega(0)=0$. Let $S\subset\R^2$ be a compact
set and $\Phi\in C(\R^2)\cap C^\infty(\R^2\setminus S)$. Suppose
that:
\begin{in.enumerate}\greek
\item\label{out.monot}
$\Phi$ is decreasing in $x$ and increasing in $y$;
\item\label{out.p>2p}
$-\Phi_x\ge 4\Phi_y>0$ on $\R^2\setminus S$;
\item\label{out.p>4omega}
$\Phi_y>4\omega'(-2\Phi_x/\Phi_y)$ on $\R^2\setminus S$;
\item\label{out.inf1} $\lim_{0<\dist((x,y),S)\to 0} \Phi_x/\Phi_y=-\infty$;
\item\label{out.null} for all \(a\leq b\) and each non-decreasing $u\in\AC[a,b]$
the sets
  $$\{x: U(x)\in S\}\text{ and }\{\Phi(U(x)): U(x)\in S\}$$ are Lebesgue null.
\end{in.enumerate}
Then there is a Lagrangian $L$ satisfying~$(\star)$ that has the
following properties:
\begin{in.enumerate}\alph
\item\label{out.2}
If $u\in\AC[a,b]$ for some \(a\leq b\in\R\), then
$$\int_a^b L(x,u,u')\, dx\ge \Phi(U(b))-\Phi(U(a)).$$
\item\label{out.3}
Equality holds in \iref{out.2} if and only if
\begin{equation}\label{u'=psi}
2\Phi_x(x,u(x))+\Phi_y(x,u(x))u'(x)=0 \textrm{ for almost every~$x\in[a,b]$}.
\end{equation}
\item\label{out.5}
Every $u\in\AC[a,b]$ satisfying \eqref{u'=psi}
is a minimizer for \(L\) on \([a,b]\).
\item\label{out.6}
If through each $(x_0,y_0)\in S$ there passes a locally absolutely
continuous solution $u\colon \R\to\R$ of~\eqref{u'=psi},
then $S$ is contained in the universal singular
set of $L$.
\end{in.enumerate}
\end{lem}

\begin{proof}
Define  auxiliary functions $\psi,\theta,\xi\in
C^\infty(\R^2\setminus S)$ by
$$\psi=-2\Phi_x/\Phi_y,\; \theta= \Phi_y-\omega'(\psi)\mbox{ and }
\xi=(-\Phi_x+\omega(\psi)-\omega'(\psi)\psi)/\theta ;$$ so that
$$\omega(\psi)+(p-\psi)\omega'(\psi)+(p-\xi)\theta= \Phi_x+p \Phi_y.$$
Note that~\eqref{out.p>2p},~\eqref{out.p>4omega} and the properties
of \(\omega\) guarantee that
  \[\theta >3\omega'(\psi)\geq 3\omega'(8)>0\text{ on }
       \R^2\setminus S,\]
 and so \(\xi\) is well-defined.

More precisely, we can use~\iref{out.p>2p} and~\iref{out.p>4omega}
to find that
  $$\Phi_y\ge\theta\ge (1-\tfrac{1}{4})\Phi_y$$
and, since \(\omega (p)-\omega'(p)p\leq \omega (0)=0\) for \(p\geq
0\),
\begin{align*}
-\Phi_x\ge -\Phi_x+\omega(\psi)-\omega'(\psi)\psi =\xi\theta &
  \ge -\Phi_x -\omega'(\psi)\psi \\
    &=(\tfrac{1}{2}\Phi_y-\omega'(\psi))\psi\\
    &\ge (\tfrac{1}{2}-\tfrac{1}{4})\psi \Phi_y=-\tfrac{1}{2}\Phi_x.
\end{align*}
Hence
\begin{equation}\label{out.ineq}
\xi\ge -\frac{1}{2}\frac{\Phi_x}{\Phi_y}\geq 2\geq 1\textrm{ and
}\psi\ge \tfrac{3}{2}\xi \geq \xi+\tfrac{1}{2}\xi \geq \xi+1\textrm{
on }\R^2\setminus S
\end{equation}
and so, by \iref{out.inf1},
\begin{equation}\label{out.inf}
\lim_{0<\dist((x,y),S)\to 0} \xi(x,y)=\infty.
\end{equation}

Let \(\gamma\) be the function given by Lemma~\ref{corner}, and
define \(F\colon\R^3\to\R\) by
$$F(x,y,p)=\begin{cases}
              \gamma(p,\xi(x,y),\theta(x,y))&\mbox{for }(x,y)\in\R^2\setminus S,\\
              0&\mbox{otherwise}.
            \end{cases}$$
Clearly, $F\in C^\infty(\R^3\setminus (S\times\R))$.
By~\eqref{out.inf}, for each $p_0\in\R$ there is an open set
$\Omega\supset S$ so that $\xi\ge p_0+1$ on $\Omega$; hence $F=0$ on
$\Omega\times (-\infty,p_0)$ and we see that $F\in C^\infty(\R^3)$.
Defining $L(x,y,p)=\omega(p)+F(x,y,p)$, it is easy to check that
\(L\) satisfies~$(\star)$.

We first note a few basic properties of $L$. Since $\xi\ge 1$,
$\omega(p)\ge\omega(0)$ and $\omega$ is strictly convex, it follows
that $L(x,y,p)\ge L(x,y,0)=\omega(0)$ and the inequality is strict
when $p\ne 0$. If $(x,y)\notin S$, we use (\ref{corner.4}), strict
convexity of $\omega$ and the definitions of $\theta$ and $\xi$, to
infer that
$$L(x,y,p)\ge \omega(\psi)+\omega'(\psi)(p-\psi)+\theta(p-\xi)=
\Phi_x+p \Phi_y$$ with equality holding if and only if $p=\psi$.

We also note some simple consequences of the assumption
\iref{out.null}. If $u\in\AC[a,b]$ is non-decreasing, then
$U(x)\notin S$ for almost every \(x\), and so $\Phi\circ U$ is
differentiable for  almost every \(x\in[a,b]\) and $(\Phi\circ
U)'=\Phi_x+\Phi_yu'$. Combining this with the properties of $L$
shown above, we see that for almost every~$x$,
$L(x,u,u')\ge(\Phi\circ U)'$ with equality holding if and only if
$u'=\psi(x,u)$. We also note that $\Phi\circ U$ has the Lusin
property --- it maps null sets to null sets: subsets of $U^{-1}(S)$
are mapped to null sets because of \iref{out.null} and null subsets
of its complement are mapped to null sets because on this set
$\Phi\circ U$ is locally absolutely continuous.

We now show that \(L\) satisfies~\iref{out.2}--\iref{out.6}.

Condition~\iref{out.monot} implies that for \(u\in \AC[a,b]\), if
$\Phi(U(a))\le \Phi(U(b))$, then $u(a)<u(b)$. Hence if such a \(u\)
is also not non-decreasing, then there is a non-decreasing
$v\in\AC[a,b]$ for which $v(a)=u(a)$, $v(b)=u(b)$ and, for almost
every~$x$, either $v'(x)=0$, or $v(x)=u(x)$ and $v'(x)=u'(x)$. Since
$\{x:0=v'(x)\ne u'(x)\}$ must have positive measure, we see that
$$\int_{a}^{b}L(x,u,u')\, dx>\int_{a}^{b}L(x,v,v')\, dx.$$
Since $L\geq 0$, it follows that, to prove~\iref{out.2} and~\iref{out.3},
we may restrict ourselves to non-decreasing $u\in\AC[a,b]$
satisfying $\Phi(U(a))\le \Phi(U(b))$. For such a \(u\), let
$G=(a,b)\setminus U^{-1}(S)$
and let $(a_j,b_j)\subset[a,b]$
be a sequence of those components of
$G$ for which $\Phi(U(a_j)) < \Phi(U(b_j))$.
Since \iref{out.null} implies that $(\Phi\circ U)(U^{-1}(S))$ is a null set,
\begin{align*}
\int_{a}^{b}L(x,u,u')\, dx
&\ge \sum_j \int_{a_j}^{b_j}L(x,u,u')\, dx\\
&\ge \sum_j \int_{a_j}^{b_j}\max(0,(\Phi\circ U)')\, dx\\
&\ge \sum_j \Phi(U(b_j))-\Phi(U(a_j))\\
&\ge \Phi(U(b))-\Phi(U(a)).
\end{align*}

The first inequality is an equality only if $G=\bigcup_j(a_j,b_j)$
and the second only if $L(x,u,u')=(\Phi\circ U)'$, and we have
already shown that this can happen only when $u'=\psi(x,u)$ almost
everywhere.

Conversely, if $u'(x)=\psi(x,u)$ almost everywhere, then $u$ is
increasing and $(\Phi\circ U)'\ge 0$ almost everywhere. Together
with the fact that \(\Phi\circ U\) has the Lusin property, this
implies that $\Phi\circ U$ is absolutely continuous,
see~\cite[Chapter~IX,\S 7.7]{saks}. Moreover, $L(x,u,u')=(\Phi\circ
U)'$ almost everywhere, and so we find $\int_{a}^{b}L(x,u,u')\,
dx=\Phi(U(b))-\Phi(U(a))$.

Statements~\iref{out.5} and~\iref{out.6} follow directly
from~\iref{out.2} and~\iref{out.3}.
\end{proof}

\begin{rem}
A slightly more general version of this lemma is obtained by
introducing a different function $\psi\in C^\infty(\R^2\setminus S)$
satisfying $\psi\ge-\Phi_x/\Phi_y$, $\Phi_y>4 \omega'(\psi)$ and
such that~\eqref{out.ineq} and~\eqref{out.inf} hold. Then there is a
Lagrangian satisfying $(\star)$ for which~\iref{out.2}--\iref{out.6}
hold with~\eqref{u'=psi} replaced by $u'(x)=\psi(x,u(x))$. This can
be used to give examples with highly non-unique minimizers.
\end{rem}

\subsection[{constructing a singular set containing a given
unrectifiable set}]{Proof of Theorem~\ref{thm.pu}: constructing a
singular set containing a given unrectifiable set}

Recall that \(\|\cdot\|_\infty\) either denotes the \(\sup\)-norm on
\(\R^2\) and \(\|f\|\) denotes the \(\sup\)-norm of \(f\), provided
\(f\) is bounded. In this section \(\lambda\) denotes Lebesgue
measure on \(\R\), and for \(a,b\in\R^2\), \([a,b]\) denotes the
closed line segment joining the points.

\begin{lem}\label{PU}
Let $S\subset\R^2$ be a compact purely unrectifiable set, $e\in\R^2$
and $\tau>0$. Then there is $f\in C^\infty(\R^2)$ for which:
\begin{itemize}
  \item $0\le f(x)\le\tau$ for all $x\in\R^2$;
  \item $\dist(\nabla f(x),[0,e])<\tau$ for all $x\in\R^2$;
  \item $\sup_{x\in S}\|\nabla f(x)-e\|_\infty < \tau$.
\end{itemize}
\end{lem}

\begin{proof}
The case $e=0$ is trivial. We can assume without loss of generality
that $e$ is the unit vector in the positive direction of the $x$
axis. Let $0<\eps<\tau$. It is enough to show that there is a
Lipschitz function $g\colon\R^2\to\R$ such that
\begin{itemize}
\item $0\le g\le \eps$ on $\R^2$;
\item $g_x\in[0,1]$ and $g_y\in[-\eps,\eps]$ at almost every point of $\R^2$;
\item $g_x=1$ in a neighbourhood of $S$.
\end{itemize}
For then a suitable mollification of $g$ gives $f$.

Let $C=1/\eps$, and let $\Omega$ be a nonempty open set containing
$S$ such that if \(\gamma\colon\R\to\R\) is a Lipschitz function
with Lipschitz constant at most \(C\), then the length of \(\{s:
(s,\gamma(s))\in \Omega\}\) is at most $\eps$. (To see that such an
open set exists, suppose instead that for each \(n\in\N\), there is
a Lipschitz function \(\gamma_n\colon \R\to\R\) for which the length
of \(P_n=\{x\in\R: (x,\gamma_n (x))\in B(S,\tfrac{1}{n})\}\) is at
least \(\eps\). On choosing a pointwise convergent subsequence of \(
(\gamma_n)\), converging to \(\gamma\), say, we find that for each
\(s\in\limsup_n P_n\), \((s,\gamma(s))\in K\), and \(\lambda
(\limsup_n P_n)>0\), contradicting the unrectifiability of \(S\),
since \(\gamma\) is Lipschitz.)

Define $g\colon\R^2\to\R$ by
$$g(x,y)=\sup\left\{x-b+\int_{(s,\gamma(s))\in \Omega}(1-
\tfrac{1}
{C}|\gamma'(s)|)\,ds\right\},$$ where the supremum is taken over all
$b\in\R$ and $\gamma\colon (-\infty,b]\to\R$, with $b\ge x$,
$\Lip(\gamma)\le C$ and $\gamma(b)=y$.

The choice $b=x$ and $\gamma (s)=y$ for \(s\leq b\) shows that
$g(x,y)\ge 0$, and it is also clear that $g(x,y)\le\sup_\gamma
\lambda\{s: (s,\gamma(s))\in \Omega\}\le\eps$. The definition of
\(g\) implies that  for arbitrary $(x,y)$ and $t>0$,
  $$g(x+t,y)\le g(x,y)+t,$$
with equality when the horizontal line segment joining $(x,y)$ and
$(x+t,y)$ lies in $\Omega$.

If $(x,y)\in\R^2$, \(b\geq x\), $t_1\ge 0$ and $|t_2|\le Ct_1$, then
the extension of the curve $\gamma\colon (-\infty ,b]\to\R$ to the
interval $[b,b+t_1]$ by a linear segment for which
$\gamma(b+t_1)=\gamma(b)+t_2$ shows that $g(x+t_1,y+t_2)\ge g(x,y)$.
In particular, for every $(x,y)$ and $t>0$, $g(x,y)\le g(x+t,y)$ and
$$g(x,y\pm t)\in
[g(x-\tfrac{1}{C}t,y),g(x+\tfrac{1}{C}t,y)]\subset[g(x,y)-\tfrac{1}{C}t,g(x,y)
+\tfrac{1}{C}t].$$
 That is, $g(x,y\pm t)-g(x,y)\in[-\eps t,\eps t]$, and so \(g\) has
 the required properties.
\end{proof}

\begin{lem}\label{P&Q}
Let $S\subset\R^2$ be a compact purely unrectifiable set, $\Omega$
an open set that contains $S$, $h^0\in C^\infty(\R^2)$,
$e^0,e^1\in\R^2$ and $\eps>0$. Then there is $h^1\in C^\infty(\R^2)$
so that:
\begin{itemize}
  \item $\|h^1-h^0\| <\eps$;
  \item $h^1=h^0$ outside $\Omega$;
  \item $\dist(\nabla h^1 (x),[e^0,e^1])
  <\eps +\|\nabla h^0(x) - e^0\|_\infty$ for \(x\in\R^2\);
and
  \item $\|\nabla h^1(x)- e^1 \|_\infty<\eps+
  \|\nabla h^0(x) - e^0\|_\infty$ for \(x\in S\).
\end{itemize}
\end{lem}

\begin{proof}
Choose $\delta>0$ so that $B(S,2\delta)\subset \Omega$. Let $\tau\in
(0,\eps/(1+\delta^{-1}))$ and take $f$ to be the function given by
Lemma~\ref{PU} for $e=e^1-e^0$. Choose $g\in C^\infty(\R^2)$ so that
  $$0\le g\le 1,\, g=1\mbox{ on }S,\, g=0\mbox{ outside }\Omega
  \mbox{ and }\|\nabla g (x)\|_\infty \le 1/\delta\text{ for }x\in\R^2.$$
Set $h^1=h^0+fg$. Clearly $h^1=h^0$ outside $\Omega$  and
  $$\sup_{x}\|h^1(x)-h^0(x)\|_\infty\le
    \sup_{x}\|f(x)\|_\infty\leq \tau <\eps.$$
Moreover for $x\in\R^2$,
\begin{align*}
\dist(\nabla h^1 (x),[e^0,e^1])&= \dist ( (\nabla
    h^0(x)-e^0)+(g\nabla f)(x) +(f\nabla g)(x), [0,e])
 \\
&\le \|\nabla h^0 (x) - e^0\|_\infty +
\dist( (g\nabla f)(x),[0,e])+  \| (f\nabla g)(x)\|_\infty\\
&<\|\nabla h^0(x) - e^0\|_\infty + \tau + \tau/\delta. \\
      &<\|\nabla
            h^0(x) - e^0\|_\infty+\eps.
\end{align*}

If $x\in S$, we replace $[e^0,e^1]$ by $e^1$ and $[0,e]$ by $e$ in
these inequalities and use that $g=1$ on \(S\) to obtain
$$\|\nabla h^1(x)- e^1 \|_\infty<\|\nabla h^0(x) -
e^0\|_\infty+\tau (1+1/\delta).$$
\end{proof}

\begin{pfof}{Theorem \ref{thm.pu}}
For \(k\geq 0\), we let
$$B_k=4+ 4\omega'(5\cdot 2^{k+4}),\, A_k=3\cdot 2^{k+2}B_k\mbox{ and }\eta_k=1-2^{-k-1},$$
and use Lemma~\ref{P&Q} recursively with vectors
$e_k=(-A_k,B_k)$ to define functions
$\Phi^k\in C^\infty(\R^2)$, open sets
$\Omega_k\subset\R^2$ and numbers $\eps_k>0$ such that:
\begin{gather}
    \label{P-0}
    \Phi^0(x,y)=-A_0 x+B_0 y,\; \Omega_0=\R^2\text{ and }\eps_0=1/4;\\
    \label{P-in}
    \|\nabla \Phi^k (x) - e_k\|_\infty<\eta_k\text{ for
    }x\in\overline{\Omega_k};\\
  \intertext{if \(a<b\), $u\in C[a,b]$ is
  non-decreasing and $\Phi\in C(\R^2)$ satisfies $\|\Phi-\Phi^k\|<2\eps_k$,
 then}
 \lambda \left(\{\Phi(x,u(x)): (x,u(x))\in \Omega_k\}\right)\leq
 1/k;\label{P-pu}\\
    \intertext{and, for \(k\ge 1\),}
    \label{P-size}
    \|\Phi^k-\Phi^{k-1}\|<\eps_{k-1};\\
    \label{P-out}
    \Phi^k=\Phi^{k-1}\text{ outside }\Omega_{k-1};\\
    \label{P-between}
    \dist(\nabla \Phi^k(x),[e_{k-1},e_k])<\eta_k\text{ for }
      x\in\overline{\Omega_{k-1}};\\
    \label{P-G} S\subset \Omega_k,\;
    \overline{\Omega_{k}}\subset B(S,2^{-k})\cap \Omega_{k-1}\text{ and }
      \eps_{k}<\eps_{k-1}/2.
\end{gather}
(We interpret \(1/0\) as \(\infty\) in~\eqref{P-pu}.)

To see that this possible, we use~\eqref{P-0} to define \(\Phi^0\),
\(\Omega_0\) and \(\eps_0\), and observe that~\eqref{P-in}
and~\eqref{P-pu} trivially hold for \(k=0\). We also observe that
$\nabla\Phi^0=e_0$, which can be considered as the appropriate
version of~\eqref{P-between} for this case. For \(k\geq 1\), define
$\Phi^k$ as the function obtained by using Lemma~\ref{P&Q} with
$\Omega=\Omega_{k-1}$, $h^0=\Phi^{k-1}$, $e^0=e_{k-1}$, $e^1=e_k$
and $\eps=\eps_{k-1}$. Then~\eqref{P-size} and~\eqref{P-out} are
just properties of \(\Phi^k\) given by the lemma.

Inequality~\eqref{P-between} follows by induction, using the
properties of \(\Phi^{k-1}\) and~\eqref{P-in} for \(k-1\),
\begin{align*}
  \dist(\nabla \Phi^k (x),[e_{k-1},e_k])
  &< \eps_{k-1}+\|\nabla\Phi^{k-1}(x)-e_{k-1}\|_\infty\\
  &<\eps_{k-1}+\eta_{k-1}\le \eta_k\text{ for }x\in\overline{\Omega_{k-1}}.
\end{align*}

For $x\in S$, we have
\begin{align*}
  \|\nabla\Phi^k (x) - e_k\| &< \eps_{k-1}+\|\nabla\Phi^{k-1} (x)-e_{k-1}\|\\
  &<\eps_{k-1}+\eta_{k-1}\le \eta_k ;
\end{align*}
by continuity these inequalities hold for some open set
$\Omega\supset S$ and so the only requirement on \(\Omega_k\) to
ensure the validity of~\eqref{P-in} is to choose \(\Omega_k\) so
that $\overline{\Omega_k}\subset\Omega$.

Since $S$ is purely unrectifiable, we may find $\delta>0$
sufficiently small so that both
$$B(S,3\delta)\subset \Omega\cap B(S,2^{-k})\cap \Omega_{k-1}$$
and for each non-decreasing $u\in C(\R)$, the linear measure of
$U(\R)\cap B(S,3\delta)$ does not exceed $1/(k(A_k+B_k+6))$. In
particular, for \(a<b\) and any non-decreasing $u\in C[a,b]$, there
are
$$a\le a_1<b_1\le a_2<b_2\le \dots\le a_q<b_q\le b$$
so that the length of each $U([a_j,b_j])$ is $\delta$ and,
$$U([a,b])\cap B(S,\delta)\subset \bigcup_{j=1}^q U([a_j,b_j])
  \subset U([a,b])\cap B(S,3\delta).$$
The second inclusion implies $q\delta\le ((A_k+B_k+6)k)^{-1}$.  In
particular, if we set $\Omega_k=B(S,\delta)$, then for
non-decreasing \(u\in C(\R)\),
\begin{equation}
  \lambda\left(\{(x,u(x)):x\in\R \}\cap
  \Omega_k\right)\leq \frac{1}{k(A_k+B_k+6)}.\label{eqn-len}
\end{equation}
Let $\eps_k=\min(\tfrac{1}{2}\eps_{k-1},\delta)$. If $\Phi$ is as
in~\eqref{P-pu}, we use the fact that the estimates
from~\eqref{P-in} are valid on $B(S,3\delta)$ to estimate
\begin{align*}
  \lambda\left(\{\Phi(x,u(x)): (x,u(x))\in \Omega_k\}\right) &=
  \lambda\left(\Phi(U([a,b])\cap B(S,\delta))\right)\\
  &\leq \sum_{j=1}^q \lambda (\Phi (U(a_j,b_j)))\\
  &\leq q(4\eps_k+(A_k+B_k+2)\delta)\\
  &\leq q(A_k+B_k+6)\delta\le 1/k.
\end{align*}
Hence~\eqref{P-pu} and~\eqref{P-size}--\eqref{P-G} hold.

By~\eqref{P-size}, the sequence $\Phi^k$ converges uniformly to some
$\Phi\in C(\R^2)$. We now show that \(\Phi\) satisfies the
hypotheses of Lemma~\ref{out}. From~\iref{P-G} and~\iref{P-out}, we
infer that $\Phi\in C^\infty(\R^2\setminus S)$. Suppose
$(x,y)\in\R^2\setminus S$, then there is $k\in\N$ such that
$(x,y)\in \overline{\Omega_{k-1}}\setminus\overline{\Omega_k}$.
By~\iref{P-out},
$$\Phi_x(x,y)=\Phi^k_x(x,y)\mbox{ and }\Phi_y(x,y)=\Phi^k_y(x,y).$$
Hence~\iref{P-between}
implies that
$$\Phi_y(x,y)\ge B_{k-1}-1\ge B_0-1\geq 3\mbox{ and }\Phi_x(x,y)\le -A_{k-1}+1\le -47;$$
in particular~\iref{out.monot} and the second inequality
of~\iref{out.p>2p}, namely \(\Phi_y>0\), hold. By~\eqref{P-between}
there is $0\le s\le 1$ so that
$$|\nabla\Phi (x,y)-(se_{k-1}+(1-s)e_k)|<1.$$
Hence
\begin{align*}
  -\Phi_x(x,y)&\le sA_{k-1}+(1-s)A_k+1
  \leq 3\cdot 2^{k+2}(sB_{k-1}+(1-s)B_k)+1\\
  &\le 3\cdot 2^{k+2}(\Phi_y(x,y)+1)+1 \le 5\cdot 2^{k+2}\Phi_y(x,y)
\end{align*}
and so \(-2\Phi_x (x,y)/\Phi_y (x,y) \leq 5\cdot 2^{k+3}\) on
\(\overline{\Omega_{k-1}}\setminus \overline{\Omega_k}\). Similarly,
\begin{align*}
-\Phi_x(x,y)&\ge sA_{k-1}+(1-s)A_k-1
\ge 3\cdot 2^{k+1}(sB_{k-1}+(1-s)B_k)-1\\
&\ge 3\cdot 2^{k+1}(\Phi_y(x,y)-1)-1 \ge 2^{k+1}\Phi_y(x,y).
\end{align*}
This gives~\iref{out.inf1} and the first inequality
of~\iref{out.p>2p}. It also shows that
$$\Phi_y (x,y)\ge B_{k-1}-1>4\omega'(5\cdot 2^{k+3})\ge
4\omega'(-2\Phi_x/\Phi_y),$$ and so~\iref{out.p>4omega} holds.

It remains to verify~\iref{out.null}. So suppose \(a<b\) and let
$u\in\AC[a,b]$ be non-decreasing. The set $\{x: (x,u(x))\in S\}$ is
null since $S$ is purely unrectifiable. Moreover,~\iref{P-size} and
$\eps_k\le\eps_{k-1}/2$ imply that for each \(k\in\N\),
$|\Phi-\Phi^k|<2\eps_k$, hence~\eqref{P-pu} implies that the measure
of $\{\Phi(x,u(x)): (x,u(x))\in S\}$ is no more than $1/k$ for all
\(k\in\N\). Hence $\{\Phi(x,u(x)): (x,u(x))\in S\}$ is also Lebesgue
null and we deduce that~\iref{out.null} holds. Hence, by
Lemma~\ref{out}, there is a Lagrangian $L$ satisfying~\( (\star)\)
and~\eqref{out.2}--\eqref{out.6}.

In order to show that $S$ is contained within the universal singular
set of \(L\), it only remains to show that the assumption
of~\eqref{out.6} holds. Let $\psi_k=-2\Phi^k_x/\Phi^k_y$ and
$\psi=-2\Phi_x/\Phi_y$. Let $(x_0,y_0)\in S$ be given. Since
$\psi_k$ is Lipschitz for each \(k\), we can find $u_k\in C^1(\R)$
so that $u_k(x_0)=y_0$ and $u_k'(x)=\psi_k(x,u_k)$. We show that
$\{u_k\}$ is an equicontinuous family. For given $\tau>0$ choose
\(q\in\N\) so that $q>2/\tau$.

Since $u_k$ is non-decreasing,~\eqref{eqn-len} implies
$$\lambda\left(\{(x,u_k(x)): (x,u_k(x))\in\Omega_q\}\right)< 1/q<\tau/2.$$
Let $\sigma= \tfrac{1}{5}2^{-(q+4)}\tau$ and consider any
$0<t-s<\sigma$. If
$$\{(x,u_k(x)): s<x<t\}\cap \Omega_q=\emptyset,$$
then, since  \(0<\psi_k\leq 5\cdot 2^{q+3}\) on
\(\R^2\setminus\Omega_q\),
$$0\le u_k(t)-u_k(s)\le 5\cdot 2^{q+3}(t-s).$$
Hence, in the general case, $0\le u_k(t)-u_k(s) <5\cdot 2^{q+3}(t-s)
+1/q<\tau$. We infer that a subsequence of $u_k$ converges locally
uniformly to a non-decreasing continuous function $u$ that satisfies
$u(x_0)=y_0$ and for which $u'=\psi(x,u)$ whenever $(x,u(x))\notin
S$. Since $S$ is purely unrectifiable, this implies that $u$ is
locally absolutely continuous, and so the hypothesis
of~\eqref{out.6} holds. Thus \(S\) is contained in the universal
singular set of \(L\).
\end{pfof}

\subsection[{a singular set meeting
 a rectifiable curve in positive length}]{Proof of Theorem~\ref{thm.non-pu}: a singular set meeting
 a rectifiable curve in positive length}

Let $A_k=4^{k+5}$, $B_k=2^{k+4}$ and recursively choose  $C_k$ to be
very large; the particular inequalities we need will follow, for
example, by setting \(C_0=0\) and picking
\begin{equation}
  C_k>8\left(\omega'(6A_{k+2}) +1+\sum_{j=0}^{k-1} C_j(1+A_j) +
(1+\sum_{j=0}^{k-1} C_j)A_k\right).\label{Ck-def}
\end{equation}
Define sets $T_k\subset\R$ and positive constants
$\ell_k$ and $\eps_k$ recursively as follows.
Let
$$T_0=[0,1],\, \ell_0=\infty\mbox{ and }\eps_0=\eps_{-1}=1.$$
For $k\ge 1$, in the $k$-th step write
$T_{k-1}$ as a finite union of
non-overlapping closed intervals $J$ each of length less than $\eps_{k-1}/A_k$
and define $T_k\subset T_{k-1}$ so that
for each of these intervals, $T_k\cap J$
is a closed interval concentric with $J$
of length $\lambda(J)(A_{k-1}-B_k)/A_k$.
Then let
$$\ell_k=\tfrac{1}{2}\min(\dist(T_k,\R\setminus T_{k-1}),\eps_{k-1})
  \text{ and }\eps_k=2^{-k-4}\ell_k/C_{k+1}.$$
Define $\chi_k\in C^\infty(\R)$ recursively by setting
$\chi_0(x)=A_0x$, and for $k\ge 1$ by defining
$$\chi_k=\chi_{k-1}\text{ outside }T_{k-1}\text{ and at the
  endpoints of the intervals $J$},$$
and by requiring
$$\chi_k'=A_k\text{ on }T_k\text{ and }B_k\le\chi_k'\le A_k\text{ on }T_{k-1};$$
the existence of $\chi_k$ is guaranteed by $0<B_k<A_{k-1}<A_k$.
Notice   that $\|\chi_k-\chi_{k-1}\|<\eps_{k-1}$.

Let $  S_k=\{(x,\chi_k(x)):x\in T_k\}$, the graph of \(\chi_k\) over
\(T_k\), and let $\Omega_0=\R^2$, and $\Omega_k=B(S_k,\ell_k)$, an
open neighbourhood around \(S_k\).

Choose functions \(\beta_k \in C^\infty(\R)\)  so that $\beta_0
\equiv 0$ and for \(k\ge 1\),
$$ \beta_k(0)=0, \beta'_k \in [0,C_k]\text{ everywhere  and }
  \beta'_k=\begin{cases}
              C_k & \text{on }(-\eps_{k-1},\eps_{k-1}),\\
              0& \text{outside }(-2\eps_{k-1},2\eps_{k-1}).
            \end{cases} $$

Choose functions $\alpha^k\in C^\infty(\R^2)$ so that
$\alpha^0\equiv \alpha^1\equiv 1$ and for  $k\geq 2$
$$ 0\le\alpha^k\le 1,\,\|\nabla\alpha^k\|\le 2/\ell_{k-1}\mbox{ and }\alpha^k=\begin{cases}
    1&\mbox{ on }\Omega_k,\\
    0&\mbox{ outside }B(\Omega_k,\ell_{k-1}).
\end{cases}$$
Define $\zeta_0(x)=\chi_0(x)=A_0x$, $\Phi^0(x,y)=0$, and for $k\ge 1$:
$$\zeta_k(x)=\chi_k(x)+\Phi^{k-1}(x,\chi_k(x))/C_k$$
and
$$\textstyle\Phi^k(x,y)=\sum_{j=0}^k \alpha^j(x,y)\beta_j(y-\zeta_j(x)).$$

An easy calculation shows that
\[\Phi_x^k (x,y)=\sum_{j=0}^k \alpha_x^j(x,y)\beta_j(y-\zeta_j(x))-
\sum_{j=0}^k \alpha^j(x,y)\beta_j'(y-\zeta_j(x))\zeta_j'(x)\]
and
\[\Phi_y^k (x,y)=\sum_{j=0}^k \alpha^j_y(x,y)\beta_j(y-\zeta_j(x))+
\sum_{j=0}^k \alpha^j(x,y)\beta_j'(y-\zeta_j(x)).\] Hence, since
$\|\nabla\alpha^0\|=\|\nabla\alpha^1\|=0$, $\|\nabla\alpha^j\| \le
2/\ell_{j-1}$ for $j\ge 2$, and $|\beta_j|\le 2C_j\eps_{j-1}$, we
estimate $\sum_{j=0}^\infty
\|\nabla\alpha^j\||\beta_j|\le\sum_{j=2}^\infty
4C_j\eps_{j-1}/\ell_{j-1}\le 1$, and so
\begin{align}
  |\Phi^k_x(x,y)+\sum_{j=0}^k \alpha^j(x,y)\beta'_j(y-\zeta_j(x))\zeta_j'(x)|
  &\le 1 \textrm{, and}\label{px1}\\
  |\Phi^k_y(x,y)-\sum_{j=0}^k \alpha^j(x,y)\beta'_j(y-\zeta_j(x))|
  &\le 1.\label{py1}
\end{align}
This implies
\begin{equation}
\|\zeta_k'-\chi_k'\|\le 1;\label{z'-c'}
\end{equation}
indeed,~\eqref{z'-c'} certainly holds for $k=0$ and hence for $k\ge
1$, we use~\eqref{px1},~\eqref{py1}, induction and~\eqref{Ck-def} to
find that for \(x\in\R\),
\begin{align*}
|\zeta_k'(x)-\chi_k'(x)|&= |(\Phi^{k-1}(x,\chi_k(x)))'|/C_k\\
&\le (|\Phi^{k-1}_x(x,\chi_k(x))|+|\Phi^{k-1}_y(x,\chi_k(x))||\chi_k'(x)|)/C_k\\
&\leq \textstyle (1+\sum_{j=0}^{k-1}\|\alpha^j\|\|\beta_j'\|\|1+\chi_j'\|
+(1+\sum_{j=0}^{k-1}\|\alpha^j\| \|\beta_j'\|)A_k)/C_k\\
&\le \textstyle
(1+\sum_{j=0}^{k-1} C_j(1+A_j) +(1+\sum_{j=0}^{k-1} C_j)A_k)/C_k
\le 1.
\end{align*}
We also observe that
\begin{align}
  |\zeta_k(x)-\chi_k(x)|  &\le\tfrac{1}{2}\eps_{k-1}
                            \text{ for }x\in T_{k-1}\text{, and}\label{z-c}\\
  \Phi^k(x,\chi_k(x)) &=0\text{ for }x\in T_{k};\label{p=0}
\end{align}
these inequalities certainly hold for $k=0$. If~\eqref{p=0} holds
for $k-1$, then for $x\in T_{k-1}$, by estimating \(\Phi_y^{k-1}\)
from~\eqref{py1} and using~\eqref{Ck-def}, we find
\begin{align*}
    |\zeta_k(x)-\chi_k(x)|&=|\Phi^{k-1}(x,\chi_k(x))-\Phi^{k-1}(x,\chi_{k-1}(x))|
                                /C_k\\
        &\le \left(1+\sum_{j=0}^{k-1} C_j\right)|\chi_k(x)-\chi_{k-1}(x)|/C_k\\
    &\le \tfrac{1}{8} C_{k}\eps_{k-1}/C_k\le \tfrac{1}{2}\eps_{k-1},
\end{align*}
which is \eqref{z-c}. Hence, since  for $x\in T_k$,
$\alpha^k(x,\chi_k(x))=1$ and
$\beta_k(\chi_k(x)-\zeta_k(x))=C_k(\chi_k(x)-\zeta_k(x))$, we have
\begin{align*}
  \Phi^k(x,\chi_k(x)) &= \Phi^{k-1}(x,\chi_k(x))+\alpha^k(x,\chi_k(x))\beta_k(\chi_k(x)-\zeta_k(x))\\
&= C_k (\zeta_k(x)-\chi_k(x))+C_k(\chi_k(x)-\zeta_k(x))=0,
\end{align*}
which is~\eqref{p=0}.

Next we note that,
since $\eps_k+\ell_{k+1}\le\ell_k $
and $\eps_{k-1}+\eps_k+\ell_{k+1}+\ell_k\le\ell_{k-1} $,
\begin{equation}
\Omega_{k+1}\subset\Omega_k\textrm{ and }
B(\Omega_{k+1},\ell_{k})\subset\Omega_{k-1}.\label{Omega1}
\end{equation}
and, since $\ell_k+\tfrac{1}{2}\eps_{k-1}\le \eps_{k-1}$, that
\begin{equation}
|y-\zeta_k(x)|\le \eps_{k-1} \textrm{ for $(x,y)\in\Omega_k.$}\label{Omega2}
\end{equation}

To deal with $(x,y)\in \Omega_q\setminus \Omega_{q+1}$
we need finer estimates of the partial derivatives
of $\Phi^k$.
Let $r:=\min(q,k)$  and $s:=\min(q+2,k)$ and observe that~\iref{Omega1} implies $\alpha^j(x,y)=0$ for
$j>q+2$, so the sums in equations~\eqref{px1} and~\eqref{py1}
finish with $j=s$. Also,
$$\left|\sum_{j=0}^{r-1} \alpha^j(x,y)\beta'_j(y-\zeta_j(x))\right|
\leq \sum_{j=0}^{r-1}\|\beta_j'\|\le \sum_{j=0}^{r-1} C_j<-1+\tfrac{1}{8}C_r$$
and, using~\iref{z'-c'},
\begin{align*}
  \left|\sum_{j=0}^{r-1} \alpha^j(x,y)\beta'_j(y-\zeta_j(x))\zeta_j'(x)\right|
  &\le \sum_{j=0}^{r-1}\|\beta_j'\| (1+\|\chi_j'\|)\\
  &\le \sum_{j=0}^{r-1} (1+A_j)C_j<-1+\tfrac{1}{8}C_r.
\end{align*}
Hence
\begin{align}\textstyle
|\Phi^k_x(x,y)+\sum_{j=r}^{s} \alpha^j(x,y)\beta'_j(y-\zeta_j(x))\zeta_j'(x)|
&\le \tfrac{1}{8}C_r\textrm{, and}\label{px2}\\
\textstyle |\Phi^k_y(x,y)-\sum_{j=r}^{s}
\alpha^j(x,y)\beta'_j(y-\zeta_j(x))| &\le
\tfrac{1}{8}C_r.\label{py2}
\end{align}
Further, using that $\alpha^r(x,y)=1$ and deducing
from~\iref{Omega2} that  $\beta'_j(y-\zeta_j(x))=C_j$ and from
\iref{z'-c'} that $\zeta_j'(x)\ge 0$ for all $j$ and $ B_r/2\le
\zeta_j'(x)\le 2A_{r+2}$ for $r\le j\le s$, we get from~\eqref{px2}
and~\eqref{py2} that
\begin{equation}
  7A_{r+2}C_{r+2}\ge -\Phi^k_x(x,y)\ge
  \tfrac{1}{4}B_rC_r\quad\text{and}\quad 4C_{r+2}\ge \Phi^k_y(x,y)\ge
  \tfrac{7}{8}C_r.  \label{psi1}
\end{equation}
Using~\eqref{py2} and~\eqref{psi1} we see that
\[\textstyle \Phi_y^k(x,y)\geq 7|\Phi^k_y(x,y)-\sum_{j=r}^{s} \alpha^j(x,y)\beta'_j(y-\zeta_j(x))|,\]
and so
\[\tfrac{6}{7}\Phi_y^k(x,y)\leq\sum_{j=r}^{s} \alpha^j(x,y)\beta'_j(y-\zeta_j(x))\leq\tfrac{8}{7}\Phi_y^k(x,y).\]
Thus
\[\tfrac{3}{7}B_r\Phi_y^k(x,y)\leq\sum_{j=r}^{s} \alpha^j(x,y)\beta'_j(y-\zeta_j(x))\zeta_j'(x)\leq\tfrac{16}{7}A_{r+2}\Phi_y^k(x,y).\]
However~\eqref{px2} and~\eqref{psi1} imply that
\[|-\Phi^k_x(x,y)-\sum_{j=r}^{s}
  \alpha^j(x,y)\beta'_j(y-\zeta_j(x))\zeta_j'(x)| \le
  \tfrac{1}{8}C_r\leq \tfrac{1}{7}\Phi_y^k(x,y),\]
and so
  \[\tfrac{1}{7}\left(-1+3B_r\right)\Phi_y^k(x,y)\leq -\Phi^k_x(x,y)
  \leq \tfrac{1}{7}\left(1+16A_{r+2}\right)\Phi_y^k(x,y),\]
giving
\begin{equation}
\tfrac{2}{7}B_r\Phi_y^k(x,y)\leq -\Phi^k_x(x,y)\leq
3A_{r+2}\Phi_y^k(x,y).\label{psi2}
\end{equation}

We are now in a position to use Lemma~\ref{out}.
Since
$$\sum_{k=1}^\infty\|\alpha^k\|\|\beta_k\|
    \leq 2\sum_{k=1}^\infty C_k\eps_{k-1}<\infty,$$
the function
$$\Phi(x,y):=\lim_{k\to\infty} \Phi^k(x,y)=
\sum_{k=0}^\infty \alpha^k(x)\beta_k(y-\zeta_k(x))$$ is continuous
on $\R^2$. Since $\sum_{k=0}^\infty \eps_k<\infty$, $\chi
:=\lim_{k\to\infty} \chi_k$ is also a continuous function on $\R$.
Moreover $\chi$ is uniformly continuous on \(\R\), since \(\chi
(x)-A_0(x)\) is zero outside \([0,1]\), and $\chi$ is  increasing,
since for all $k\geq 1$ and $y>x$
$$\chi_k(y)-\chi_k(x)\ge\min (A_0,B_1)(y-x)= B_1 (y-x).$$

Let $T=\bigcap_{k=0}^\infty T_k$ and $S=\{(x,\chi(x)): x\in T\}$.
Then $\lambda(T)=0$, since
$\lambda(T_k)=(A_{k-1}-B_k)\lambda(T_{k-1})/A_k$ and
$\prod_k(A_{k-1}-B_k)/A_k=0$. Similarly, we observe that
$$\lambda(\chi_k(T_k))=A_k\lambda(T_k)=(A_{k-1}-B_k)\lambda(T_{k-1})=
(A_{k-1}-B_k)\lambda(\chi_{k-1}(T))/A_{k-1},$$
$\prod_k(A_{k-1}-B_k)/ A_{k-1}>0$, and, since \(\chi_{k+1}
(T_{k+1})\subset \chi_k (T_k)\), \(\chi (T)=\cap_{k=1}^\infty \chi_k
(T_k)\). Hence $\lambda(\chi(T))>0$, and it follows that $S$ is a
rectifiable set of positive linear measure. Note also
that~\iref{Omega1} implies that $S=\bigcap_{k=0}^\infty \Omega_k$.
In particular, the sum defining $\Phi$ is locally finite on
$\R^2\setminus S$. This means that $\Phi\in C^\infty(\R^2\setminus
S)$ and that~\iref{out.monot},~\iref{out.p>2p}, \iref{out.p>4omega}
and \iref{out.inf1} follow immediately from~\eqref{Ck-def},
\eqref{psi1} and~\eqref{psi2}. The condition~\iref{out.null} holds
because for any $u\colon \R\to\R$, $\{x: (x,u(x))\in S\}\subset T$
is a null set and $\{\Phi(x,u(x)): (x,u(x))\in S\}\subset\{0\}$,
since~\eqref{p=0} implies that $\Phi(x,y)=0$ on $S$.

Finally, to verify that the assumptions of~\iref{out.6} hold, for
$k\ge 1$ let
$$\psi_k=-2\Phi^k_x/\Phi^k_y\mbox{ and }\psi=-2\Phi_x/\Phi_y.$$
Fix $(x_0,y_0)\in S$.

Since $\psi_k$ is Lipschitz for each \(k\), we can
find $u_k\in C^1(\R)$ so that $u_k(x_0)=y_0$ and
$u_k'(x)=\psi_k(x,u_k)$.
We show that $\{u_k\}$ is an equicontinuous family.
For given $\tau>0$ choose $q$ so that
$$2\ell_q +2\sum_{j\ge q}\eps_j<\tau/6$$
and use the uniform continuity of \(\chi \) to choose
$\sigma\leq\tfrac{1}{18}\tau/A_{q+2}$ so that for \(0\leq
t-s<\sigma\),
$$0\leq\chi(t)-\chi(s)<\tau/6.$$
Consider $0<t-s<\sigma$. If
$$\{(x,u_k(x)): s<x<t\}\cap\Omega_q=\emptyset,$$
we have, since \(\psi_k\) is at most \( 6A_{q+2}\),
$$0\le u_k(t)-u_k(s)\le 6A_{q+2}(t-s)<\tau/3.$$
If $(s,u_k(s))\in\overline{\Omega_q}$
and $(t,u_k(t))\in\overline{\Omega_q}$,
we have
$$0\le u_k(t)-u_k(s)\le\chi_q(t)-\chi_q(s)+2\ell_q
\le \chi(t)-\chi(s)+2\ell_q +2\sum_{j\ge q}\eps_j<\tau/3.$$
In the general case, the interval $[t,s]$ can be written as
the union of three non-overlapping intervals each of which satisfies
one of these conditions, so $0\le u_k(t)-u_k(s)<\tau$.
We infer that a subsequence of $u_k$
converges locally uniformly to a non-decreasing continuous function $u$
that satisfies $u(x_0)=y_0$
and for which $u'=\psi(x,u)$ whenever $(x,u(x))\notin S$.
Note that
\begin{align*}
\frac{d}{dx}\Phi^k(x,u_k)&= \Phi^k_x(x,u_k(x))+\Phi^k_y (x,u_k(x))u_k'(x) \\
&=\Phi^k_x (x,u_k)+\Phi^k_y (x,u_k)\times\frac{-2\Phi^k_x (x,u_k)}{\Phi_y^k (x,u_k)}\\
&= -\Phi_x^k (x,u_k(x))\ge 1.
\end{align*}
Hence $\Phi^k(t,u_k(t))-\Phi^k(s,u_k(s))\ge t-s$, and so
$\Phi(t,u(t))-\Phi(s,u(s))\ge t-s$ for $t>s$. Since $\Phi(x,y)=0$ on
$S$, the graph of $u$ meets $S$ only at the point $(x_0,y_0)$. Hence
$u'=\psi(x,u)$ for $x\in \R\setminus\{x_0\}$ and the monotonicity of
$u$ implies that it is locally absolutely continuous.
Lemma~\ref{out} now guarantees the existence of a Lagrangian with
the required properties.

\subsection[{a continuous Lagrangian with residual universal singular
set}]{Proof of Theorem~\ref{resid}: a continuous Lagrangian with
residual universal singular set}\label{residuss}

This Lagrangian is constructed by a variant of the general
construction of Lagrangians with large universal singular sets that
we described earlier; the technicalities are simpler in one sense,
since we do not have to care about the smoothness of the Lagrangian,
however complications are caused by the fact that we have to work
with derivatives of non-smooth functions.

\begin{lem}\label{resid1.lem}
Suppose that $z_0\in\R^2$, $\eps>0$ and $0\ne P\in\R^2$. Then there
are an open set $G\subset \R^2$, a Lipschitz function
$\Phi\colon\R^2\to\R$ and a continuous function
$\phi\colon\R^2\to\R^2$ such that
\begin{in.enumerate}
\item\label{r1.0}
$z_0\in\overline{G}$, $\diam(G)<\eps$ and $\partial G$ is Lebesgue
null;
\item\label{r1.2}
$|\Phi|\le\eps$ on \(\R^2\), $\|\phi\|\le\eps$, $\nabla\Phi=\phi$ on
$\R^2\setminus\overline{G}$ and $\nabla\Phi=P$ on $G$;
\item\label{r1.1}
if $\|z-z_0\|\ge\eps$, then  $|\Phi(w)-\Phi(z)|\le\eps \|w-z\|$ for
all \(w\in\R^2\);
\item\label{r1.4}
if $\gamma\colon [a,b]\to\R^2$ is absolutely continuous, then
\begin{enumerate}
\item $\frac{d}{dx}\Phi(\gamma(x)) =P\cdot\gamma'(x)$ for a.e.\ $x$
for which $\gamma(x)\in\overline G$, and
\item $\frac{d}{dx}\Phi(\gamma(x))=\phi(\gamma(x))\cdot\gamma'(x)$
for a.e.\ $x$ for which $\gamma(x)\notin G$.
\end{enumerate}
\end{in.enumerate}
\end{lem}

\begin{proof}
Since the statement does not depend on the choice of origin and
coordinates, we may assume that $z_0=(0,0)$ and $P=(0,c)$, where
$c>0$.

Let $0<\delta<1$ be small enough; the exact conditions are
determined below. Choose $f\in C^1(\R)$ so that
 $$f(x)=0\text{ for }|x|\ge\delta,\, 0<f(x)<\delta\text{ for }x\in(-\delta,\delta)
  \text{ and }|f'(x)|<\delta\text{ for all }x.$$
Let $G=\{(x,y): 0<y<f(x)\}$, and define
$$\Phi(x,y)=\begin{cases}
0& \textrm{ if $y\le 0$},\\
cy &\textrm{ if $0<y\le f(x)$},\\
cf(x)& \textrm{ if $y > f(x)$},
\end{cases}
$$
and, for $(x,y)\in\R^2\setminus G$, let
$$\phi(x,y)=\begin{cases}
(0,0)& \textrm{ if $y\le 0$},\\
(cf'(x),0)& \textrm{ if $y\ge f(x)$.}
\end{cases}
$$

Then~\iref{r1.0} holds, provided $3\delta<\eps$. To ensure
$|\Phi|\le\eps$ on \(\R^2\),  we require $\delta c\le\eps$. The same
assumption on $\delta$ also guarantees $\|\phi\|\le\eps$ on
$\R^2\setminus G$. Since $\phi$ is (well-defined and) continuous on
$\R^2\setminus G$, we can then extend it to  a continuous function
on $\R^2$ with $\|\phi\|\le\eps$. The remainder of~\iref{r1.2} is
obvious.

If $[z,w]$ does not meet $\overline{G}$,~\iref{r1.1} holds, since
$\|\nabla\Phi (x)\|\le\eps$ for every $x\in [z,w]$. If $[z,w]$ meets
$\overline{G}$, then the intersection has length at most $3\delta$
and $\|z-w\|\ge\eps-2\delta$. Since $\|\nabla\Phi (x)\|\le c$ for
\(x\in [z,w]\cap\overline{G}\) and $\|\nabla\Phi (x)\|\le \delta c$
for \(x\not\in [z,w]\cap\overline{G}\), we use that \(\|w-z\|\geq
\|z-z_0\|-\diam(G)\geq\eps -2\delta\) to estimate that
  $$|\Phi(w)-\Phi(z)|\le 3\delta c+ \delta c \|w-z\|\le \eps\|w-z\|,$$
provided that  $3\delta c\le \tfrac{1}{2}\eps(\eps-2\delta)$ and
$\delta c\le\tfrac{1}{2}\eps$.

Finally, if $\gamma=(\gamma_1,\gamma_2)$ is as in~\iref{r1.4}, then,
since $\Phi(\gamma(x)) =P\cdot\gamma(x)$ whenever
$\gamma(x)\in\overline{G}$, we conclude that
$\frac{d}{dx}\Phi(\gamma(x)) =P\cdot\gamma'(x)$ whenever
$\gamma(x)\in\overline{G}$, $x$ is not an isolated point of
$\gamma^{-1}(\overline{G})$, and $\gamma$ and \(\Phi\circ\gamma\)
are differentiable at $x$. Similarly, we infer from $\Phi(\gamma(x))
=cf(\gamma_1(x))$ when $\gamma_2(x)\ge f(\gamma_1(x))$, that
$\frac{d}{dx}\Phi(\gamma(x))=cf'(\gamma_1(x))\gamma_1'(x)=
\phi(\gamma(x))\cdot\gamma'(x)$ for almost all $x$ for which
$\gamma_2(x)\ge f(\gamma_1(x))$, and, from $\Phi(\gamma(x)) =0$ when
$\gamma_2(x)\le 0$, that $\frac{d}{dx}\Phi(\gamma(x))=0=
\phi(\gamma(x))\cdot\gamma'(x)$ for almost all $x$ for which
$\gamma_2(x)\le 0$.
\end{proof}

In the following,   $\HH(A)$ denotes the one-dimensional Hausdorff
capacity of a set $A\subset\R^2$; that is,
$$\HH(A)=
  \inf\left\{\sum_{j=1}^\infty \diam(A_j): A\subset
  \bigcup_{j=1}^\infty A_j\right\}.$$
Recall that for \(u\colon A\subset\R\to\R\), \(U\) denotes the
function \(x\mapsto (x,u(x))\).

\begin{lem}\label{resid2.lem}
  Let $\emptyset\not=H\subset\R^2$ be open, $\eps>0$ and
$P=(-A,B)\in\R^2$ with $A,B>0$. Then there is an open set $G\subset
H$ that is dense in $H$, a Lipschitz function $\Phi\colon\R^2\to\R$,
and a continuous function $L\colon\R^3\to\R$ for which
\begin{in.enumerate}\alph
  \item\label{r2.0} $\HH(G)\le \eps$;
  \item\label{r2.1} $|\Phi|\le \eps$ on \(\R^2\);
  \item\label{r2.2} $|L(x,y,p)-\max(0,-A+Bp)|\le \eps(1+|p|)$;
  \item\label{r2.5} $\int_a^b L(x,u,u')\,dx\ge\Phi(U(b))-\Phi(U(a))$
  for every $u\in\AC[a,b]$;
  \item\label{r2.7}
  if $u\in\AC[a,b]$ satisfies $u'(x)\ge A/B+1$ for a.e.~$x$ for
  which $U(x)\in G$, and there is an absolutely equicontinuous
  sequence $u_j\in\AC[a,b]$ converging to $u$ such that $|u_j'(x)|\le
  A/B-1$ for a.e.~$x$ for which $U_j(x)\notin G$, then
  $\int_a^b L(x,u,u')\,dx=\Phi(U(b))-\Phi(U(a))$.
\end{in.enumerate}
\end{lem}

\begin{proof}
We may assume that $2\eps(2+A/B)\le \min(A,B)$. Let $z_1,z_2,\ldots
\in H$ be a sequence that is dense in $H$.

We  use Lemma~\ref{resid1.lem} recursively to define open sets
$G^k\subset\R^2$ that have pairwise disjoint closures, Lipschitz
functions $\Phi^k\colon\R^2\to\R$ and continuous functions
$\phi^k\colon\R^2\to\R^2$. In addition, we also define real numbers
\(\eps_k>0\) and continuous functions $\zeta^k\colon \R^2\to\R^2$
such that
\begin{equation}
  \sum_{j=1}^k \nabla\Phi^j=\sum_{j=1}^k \zeta^j\text{ on }
    \bigcup_{j=1}^k G^j\text{ for each }k,\label{eqn-sumphi=sumzeta}
\end{equation}
and
  $$\|\zeta^k (x)\|\le \eps_k\text{ for }x\in\R^2\text{ and }k\ge
2.$$

To start the recursion, we choose $0<\eps_1 <\tfrac{1}{2}\eps $ so
that $\{z: \|z-z_{1}\|\le\eps_1\}\subset H$, let
$P_1=(-A_1,B_1)=(-A,B)$, let $\zeta^1(z)=(-A,B)$, and use
Lemma~\ref{resid1.lem} with $z_0=z_{1}$, $\eps=\eps_1$ and $P=P_1$
to define $G^1,\Phi^1$ and $\phi^1$.

Assume now that $k\ge 2$ and that all objects are defined for $1\le
j<k$. Let $i_k$ be the first index for which
$z_{i_k}\notin\bigcup_{j=1}^{k-1} \overline{G^j}$, and choose
$0<\eps_k <2^{-k}\eps $ so that
  $$\{z:\|z-z_{i_k}\|\le\eps_k\}\subset H\setminus
    \bigcup_{j<k}\overline{G^j}.$$
Let
  $$P_k=(-A_k,B_k)=\sum_{j=1}^{k-1}(\zeta^j(z_{i_k})-\nabla\Phi^j(z_{i_k})).$$
(Notice that $\|P_k-P\|\le 2\eps$, and so $A_k,B_k>0$.) Since
\(z\mapsto \sum_{j=1}^{k-1}(\zeta^j(z)-\nabla\Phi^j(z))\) is
continuous at \(z_{i_k}\), there is \(0<\eps_k'\leq \eps_k\) so that
  \[\left\|P_k-\sum_{j=1}^{k-1}(\zeta^j(z)-\nabla\Phi^j(z))\right\|
  \leq\eps_k\text{ on }B(z_{i_k},\eps_k').\]
Now use Lemma~\ref{resid1.lem} with $z_0=z_{i_k}$, $\eps=\eps_k'$
and $P=P_k$ to define $G^k$, $\Phi^k$ and $\phi^k$. Let
  $$\zeta^k(z)=
        \begin{cases}
          \nabla\Phi^k(z)&\text{on }\bigcup_{j=1}^{k-1}
          \overline{G^j}\text{, and}\\
          P_k-\sum_{j=1}^{k-1}(\zeta^j(z)-\nabla\Phi^j(z))
          &\text{for }z\in \overline{G^k}.
        \end{cases} $$
Then $\zeta^k\colon \bigcup_{j=1}^{k} \overline{G^j}\to\R^2$ is
continuous and $\|\zeta^k\|\le\eps_k$. Hence we may extend
\(\zeta^k\) to a continuous function $\zeta^k\colon\R^2\to\R^2$ for
which $\|\zeta^k\|\le\eps_k$. It remains to
verify~\eqref{eqn-sumphi=sumzeta}: on $\bigcup_{j=1}^{k-1} G^j$ we
have
  $$\sum_{j=1}^k \zeta^j =\zeta^k+\sum_{j=1}^{k-1} \nabla\Phi^j
    =\sum_{j=1}^{k} \nabla\Phi^j,$$
 and on $G_k$,
  $$\sum_{j=1}^k
    \zeta^j=P_k+ \sum_{j=1}^{k-1} \nabla\Phi^j =\sum_{j=1}^{k}
    \nabla\Phi^j.$$

Since $|\Phi^k|\le\eps_k$ on \(\R^2\), the series
$\Phi:=\sum_{j=1}^\infty\Phi^j$ defines a continuous function on
$\R^2$. Moreover, $|\Phi|\le\eps$ (so~\iref{r2.1} holds) and $\Phi$
is Lipschitz, since the sets $G^j$ are disjoint, and so each partial
sum $\sum_{j=1}^k\Phi^j$ is Lipschitz. Also observe that
$\left\|\nabla\left(\sum_{j=1}^k\Phi^j\right)(x)\right\|\leq
A+B+2\eps$ for almost every \(x\), since
\(\|\nabla\Phi^j\|\leq\eps_j\) outside \(\overline{G^j}\),
\(\nabla\Phi^j=P_j\) on \(G^j\), \(\|P_j-P\|\leq 2\eps\) and
\(\|P\|\leq A+B\).

 Clearly the set $G=\bigcup_{j=1}^\infty G^j$ is an open dense subset of
 $H$, and $\HH(G)\le\sum_{j=1}^\infty\eps_k\le\eps$ so~\iref{r2.0} holds.

Write  $\phi^k=(\phi^k_1,\phi^k_2)$ and
$\zeta^k=(\zeta^k_1,\zeta^k_2)$, and define functions
$\LO,\LI\colon\R^3\to\R$ by
\begin{align*}
  \LO(x,y,p)&=\sum_{j=1}^\infty
  \left(\phi^j_1(x,y)+p\phi^j_2(x,y)\right)
  \intertext{ and}
  \LI(x,y,p)&=\sum_{j=1}^\infty\left(\zeta^j_1(x,y)+p\zeta^j_2(x,y)\right).
\end{align*}
Recalling that $\|\phi^k\|\le\eps_k$ for all $k$,
$\|\zeta^k\|\le\eps_k$ for $k\ge 2$ and $\zeta^1=(-A,B)$, we see
that $\LO$, $\LI$ are continuous and satisfy $|\LO|\le\eps(1+|p|)$
and $|\LI -(- A + Bp)|\le\eps(1+|p|)$ on \(\R^3\). It follows that
the function
$$L(x,y,p):=\max\left(\LI(x,y,p),\LO(x,y,p)\right)$$
satisfies~\iref{r2.2}. To prove the remaining statements, we first
show that
\begin{align}
  L(x,y,p)&=\LO(x,y,p) \text{, when }|p|\le A/B-1
  \label{r.p.lo}\\
  \intertext{and}
  L(x,y,p)&=\LI(x,y,p) \text{, when }(x,y)\in G\text{ and }p\ge
  A/B+1.
  \label{r.p.li}
\end{align}
Indeed, if $|p|\le A/B-1$, then, since $2\eps\leq B^2/A$,
  $$\LI(x,y,p)\le -A+Bp +\eps(1+|p|) \le -\eps(1+|p|)\le \LO(x,y,p),$$
and if $p\ge A/B+1$, then, since \(2\eps (2+A/B)\leq B\),
  $$\LI(x,y,p)\ge -A+Bp -\eps(1+|p|) \ge \eps(1+|p|)\ge \LO(x,y,p).$$

Let $Z_k=\bigcup_{j=k}^\infty\{z\in\R^2:\|z-z_{i_j}\|<\eps_j\}$,
$Z=\bigcap_{k=1}^\infty Z_k$ and observe that if $z\notin Z_k$ then
$|\Phi^j(w)-\Phi^j(z)|\le \eps_j\|w-z\|$ for all $w\in\R^2$ and
$j\ge k$. Since $\Phi^k$ is differentiable on $\R^2\setminus
\partial G^k$, we infer that
\begin{equation}\label{r.p.d1}
\nabla\Phi(z)=\sum_{j=1}^\infty\nabla\Phi^j(z) \text{, when }
z\notin Z\cup \bigcup_{j=1}^\infty\partial G^j.
\end{equation}
The same argument, together with the fact that $\partial G^k\cap
Z_{k+1}=\emptyset$, shows that
\begin{equation}\label{r.p.d2}
\nabla\left(\sum_{j\ne k} \Phi^j\right)(z)=\sum_{j\ne k}
\nabla\Phi^j(z) \text{, when } z\in\partial G^k.
\end{equation}

If $u\in\AC[a,b]$, then the function $\Phi\circ U$ is absolutely
continuous, since $\Phi$ is Lipschitz. Noting that $\HH(Z)=0$, and
so $U(x)\notin Z$ for almost every $x$, we see that its derivative
is described by one of the following three cases for almost every
$x$:
\begin{enumerate}
\item\label{r.p.G}
If $U(x)\in G$, then~\eqref{r.p.d1} and~\eqref{eqn-sumphi=sumzeta}
give $\frac{d}{dx}\Phi(U(x))=\LI(x,u,u').$
\item\label{r.p.C}
If $U(x)\notin Z\cup \bigcup_{j=1}^\infty\overline{G^j}$,
then~\eqref{r.p.d1} gives that $\frac{d}{dx}\Phi(U(x))=\LO(x,u,u').$
\item\label{r.p.D}
If $U(x)\in\partial\left( \bigcup_{j=1}^\infty G^j\right)$,
then~\eqref{r.p.d2} and~\iref{r1.4} give that
$\frac{d}{dx}\Phi(U(x))=\LI(x,u,u')=\LO(x,u,u').$
\end{enumerate}

The statement~\iref{r2.5} follows immediately from this and  the
definition of $L$. To deduce~\iref{r2.7}, let $u$ satisfy its
assumptions. Since $u'\ge A/B+1$ for a.e.~$x$ for which $U(x)\in
G$,~\iref{r.p.G} and~\eqref{r.p.li} imply that
$\frac{d}{dx}\Phi(x,u)=L(x,u,u')$ for such~$x$. By~\iref{r.p.D}, the
same expression for the derivative of $\Phi(x,u)$ holds for a.e.~$x$
for which $U(x)\in \partial G_k$ for some $k$. Hence it is enough to
show that it holds for a.e.~$x$ for which $U(x)\notin
\bigcup_{j=1}^\infty\overline{G^j}$, since then it will hold almost
everywhere implying that $\F(u;a,b)= \Phi(b,u(b))-\Phi(a,u(a)).$

Let $W_k$ be the $x$-projection of $\bigcup_{j=k+1}^\infty
\overline{G^j}$. If $K$ is a compact subset of $\{x: U(x)\notin
\bigcup_{j=1}^\infty\overline{G^j}\}$ and $k$ is fixed, then for all
sufficiently large $l$, $U_l(x)\notin \bigcup_{j=1}^k\overline{G^j}$
for all $x\in K$. Hence $u_l'(x)\le A/B-1$ for almost all $x\in
K\setminus W_k$. Since $u_l'$ converges to $u'$ weakly in
$L^1[a,b]$, we infer that $u'(x)\le A/B-1$ for almost all $x\in
K\setminus W_k$. Since $k$ is arbitrary and the measure of $W_k$
tends to zero, we have that $u'(x)\le A/B-1$ for almost all $x\in
K$. Finally, we use that $K$ is arbitrary to deduce that $u'(x)\le
A/B-1$ for almost all $x$ for which $U(x)\notin \bigcup_{j=1}^\infty
\overline{G^j}$ and conclude from~\iref{r.p.C} and~\eqref{r.p.lo}
that $\frac{d}{dx}\Phi(x,u(x))=L(x,u,u')$ for almost all $x$ from
this set, as required.
\end{proof}

\begin{proof}[Proof of Theorem \ref{resid}]
Let $\omega\ge 0$ be convex and superlinear. We let
  $$A_k=1+\omega(2^{2k+2}),\, B_k=2^{-2k+1}A_k\text{ and }
    \eps_k=2^{-k}(1+2^{2k+2})^{-1}.$$
To start the recursive construction, let $G_0=\R^2$,
$\Phi_0(x,y)=-A_0x+B_0y$ and $L_0(x,y,p)=-A_0+B_0p$. (For future
use, notice that $\int_a^b L_0(x,u,u')\,dx =
\Phi_0(U(b))-\Phi_0(U(a))$ for every $u\in\AC[a,b]$.) Using
Lemma~\ref{resid2.lem} recursively with $H=G_{k-1}$, $\eps=\eps_k$
and $P=(-A_k,B_k)$, we define open sets $G_k\subset G_{k-1}$ that
are dense in $G_{k-1}$, Lipschitz functions $\Phi_k\colon
\R^2\to\R$, and continuous functions $L_k\colon\R^3\to\R$ with the
properties described there.

By~\iref{r2.1}, $\Phi(x,y):=\sum_{k=0}^\infty\Phi_k(x,y)$ is
continuous  on $\R^2$. If $|p|< 2^{2j-1}$, then~\iref{r2.2} implies
that $|L_k(x,y,p)|\le\eps_k(1+2^{2j-1})$ for $k\ge j$. Hence both
\begin{align*}
\tilde L(x,y,p)&:=\sum_{k=0}^\infty \max(-2^{-k}, L_k(x,y,p))
\intertext{and}
 L(x,y,p)&:=\max(\omega(p),\tilde L(x,y,p))
\end{align*}
are continuous functions on $\R^3$.

Clearly $L$ is superlinear. Noting that
  $$\max (0,-A_j+B_j p)=\begin{cases}
    0 &\text{when }0\leq p\leq A_j/B_j,\\
    -A_j+B_jp &\text{when }p\geq A_j/B_j,
  \end{cases} $$
we deduce from~\iref{r2.2} that for \(0\leq p\leq  A_j/B_j\),
  $$ L_j(x,y,p)\ge -\eps_j(1+A_j/B_j) $$
and for \( p > A_j/B_j\),
\begin{align*}
  L_j (x,y,p) & \geq B_j p-A_j-\eps_j (1+p)\\
  &= (B_j-\eps_j)(p-A_j/B_j)-\eps_j (1+A_j/B_j)
  &\geq -\eps_j (1+A_j/B_j) .
\end{align*}
Hence for any \(p\geq 0\),
\begin{equation}
  L_j(x,y,p)\ge -\eps_j(1+A_j/B_j)\ge -2^{-j},\label{eqn-Lj}
\end{equation}
   and so for
\(p\in[2^{2k},2^{2k+2}]\), we use~\iref{r2.2} for \(L_k\), the fact
that \(L_0 (x,y,p)\geq 0\) and~\eqref{eqn-Lj} to estimate
  \begin{align*}
    \tilde L(x,y,p)&= \sum_{j=0}^\infty L_j (x,y,p)\\
    &= L_k (x,y,p) +\sum_{j=1,j\not=k}^\infty L_j (x,y,p)\\
    &\ge -A_k+pB_k
    -\sum_{j=1}^\infty\eps_j(1+2^{2j+2})\\
  &\ge A_k-1= \omega(2^{2k+2})\ge\omega(p).
  \end{align*}
Hence for \(p\geq 1\),
\begin{equation}\label{res.pf.2}
L(x,y,p)=\sum_{j=0}^\infty L_j(x,y,p) \textrm{ and } L_j(x,y,p)\ge
-2^{-j} .
\end{equation}

If $u\in\AC[a,b]$, we have $\int_a^b L_j(x,u,u')\,dx\ge
\Phi_j(U(b))-\Phi_j(U(a))$ by~\iref{r2.5}, and so
\begin{align}
  \F(u;a,b)&\ge
  \int_a^b \sum_{j=0}^\infty \max(-2^{-j},L_j(x,u,u'))\,dx\nonumber\\
  &=\sum_{j=0}^\infty\int_a^b \max(-2^{-j},L_j(x,u,u'))\,dx\nonumber\\
  &\ge \sum_{j=0}^\infty\left(\Phi_j(U(b))-\Phi_j(U(a))\right)\nonumber\\
  &=\Phi(U(b))-\Phi(U(a))\label{res.pf.1}
\end{align}

Let $G=\bigcap_{k=0}^\infty G_k$. Since $G_0=\R^2$ and $G_k$ is
dense in $G_{k-1}$, $G$ is a residual subset of $\R^2$. We show that
$G$ is contained in the universal singular set of $L$. For this,
assume that $(x_0,y_0)\in G$ and define
$\psi\colon\R^2\to[1,\infty]$ by $\psi(z)=A_k/B_k-1=2^{2k-1}-1$ for
$z\in G_{k-1}\setminus G_k$ and $\psi(z)=\infty$ for $z\in G$. Then
$\psi$ is lower semicontinuous.

Let $\psi_k\colon\R\to[1,\infty)$ be bounded continuous functions
such that $\psi_k\nearrow\psi$ as $k\to\infty$. Then the equation
$u'=\psi_k(x,u)$ has a global $C^1$ solution $u_k$ such that
$u_k(x_0)=y_0$.

Denote by $f(x)$ the supremum of those values $A_k/B_k$ for which
$(x,y)\in G_{k-1}$ for some $y$. Then
  $$\int_a^b f(x)\,dx\le
    (b-a)A_1/B_1+\sum_{k=2}^\infty \HH(G_{k-1})A_k/B_k<\infty$$
(the convergence of the series follows from~\iref{r2.0}), so $f$ is
locally integrable. By definition, $\psi(x,y)\le f(x)$ for every
$(x,y)\in \R^2$. Since $\psi_k\le\psi$, we infer that $0\le u_k'\le
f$, and conclude that the sequence $u_k$ is locally absolutely
equicontinuous. Hence it has a subsequence converging (locally
uniformly) to a locally absolutely continuous function
$u\colon\R\to\R$ such that $u(x_0)=y_0$. Since $u_l'=
\psi_l(x,u_l)\ge \psi_k(x,u_l)$ for $l\ge k$ and since
$\psi_k(x,u_l)$ converges locally uniformly to $\psi_k(x,u)$ as
$l\to\infty$, we have $u'\ge\psi_k(x,u)$ for all $k$ and so
$u'\ge\psi(x,u)$.

Since $u'(x)=\psi(x,u)\ge A_k/B_k+1$ when $(x,u)\in G_k$ and
$$0\leq u'_l(x)\leq \psi_l(x,u_l)\le \psi (x,u_l)\leq
A_k/B_k-1\text{ when }(x,u)\notin G_k,$$
 $u$ satisfies the assumptions
of~\iref{r2.7} for each interval $[a,b]$ and each $k$. Hence
$\int_a^b L_k(x,u,u')\,dx=\Phi_k(b,u(b))-\Phi_k(a,u(a))$. Since
$u'\ge\psi(x,u)\ge 1$, \eqref{res.pf.2} implies that
$L(x,u,u')=\sum_{j=0}^\infty L_j(x,u,u')$ and $L_j(x,u,u')\ge
-2^{-j}$. This justifies exchange of integration and summation,
allowing us to conclude that $\int_a^b L(x,u,u')\,dx=
\sum_{j=0}^\infty \int_a^b L_j(x,u,u')\,dx= \Phi(U(b))-\Phi(U(a))$
for any $a<b$ and so, because of~\eqref{res.pf.1},  $u$ is a
minimizer on any interval $[a,b]$. Finally, since
$$\lim_{x\to x_0} u'(x)\ge \lim_{(x,y)\to
(x_0,y_0)} \psi(x,y)=\infty,$$ we have $u'(x_0)=\infty$ and we are
done.
\end{proof}

\end{document}